\documentclass[oneside, reqno]{amsart}

\usepackage{xypic}
\input xy
\xyoption{all}
\usepackage{epsfig}
\usepackage{amsthm}
\usepackage{amssymb}
\usepackage{amsmath}
\usepackage{amscd}
\usepackage{color}
\usepackage[T1]{fontenc}
\usepackage{stackengine}
\usepackage{upgreek}
\usepackage[font=scriptsize]{caption}
\usepackage{wrapfig}
\stackMath

\usepackage[parfill]{parskip}
\usepackage{tikz,pgfplots,caption}
\usetikzlibrary{patterns}

\newcommand{\bg}{\begin{equation}}
\newcommand{\ed}{\end{equation}}
\newcommand{\bga}{\begin{eqnarray}}
\newcommand{\eda}{\end{eqnarray}}

\def\cbdu{\par{\raggedleft$\Box$\par}}

\newcommand \ea{\varepsilon^+}
\newcommand \eb{\varepsilon^-}
\newcommand \dea{\delta^+}
\newcommand \deb{\delta^-}

\newtheorem {Theorem}  {Theorem}

\numberwithin{Theorem}{section}

\newtheorem {Lemma}[Theorem]  {Lemma}

\theoremstyle{definition}
\newtheorem{Definition}[Theorem]{Definition}
\theoremstyle{remark}
\newtheorem{Remark}[Theorem]{\bf Remark}
\newtheorem{Conjecture}[Theorem]{\bf Conjecture}

\expandafter\chardef\csname pre amssym.def
at\endcsname=\the\catcode`\@ \catcode`\@=11
\def\undefine#1{\let#1\undefined}
\def\newsymbol#1#2#3#4#5{\let\next@\relax
 \ifnum#2=\@ne\let\next@\msafam@\else
 \ifnum#2=\tw@\let\next@\msbfam@\fi\fi
 \mathchardef#1="#3\next@#4#5}
\def\mathhexbox@#1#2#3{\relax
 \ifmmode\mathpalette{}{\m@th\mathchar"#1#2#3}%
 \else\leavevmode\hbox{$\m@th\mathchar"#1#2#3$}\fi}
\def\hexnumber@#1{\ifcase#1 0\or 1\or 2\or 3\or 4\or 5\or 6\or 7\or 8\or
 9\or A\or B\or C\or D\or E\or F\fi}

\font\teneufm=eufm10 \font\seveneufm=eufm7 \font\fiveeufm=eufm5
\newfam\eufmfam
\textfont\eufmfam=\teneufm \scriptfont\eufmfam=\seveneufm
\scriptscriptfont\eufmfam=\fiveeufm

\catcode`\@=\csname pre amssym.def at\endcsname

\newcounter{remark}
\newcommand{\R}{\mathbf{R}}

\def  \R   {{\mathbb R}}

\def  \12  {{\frac{1}{2}}}

\def\build#1_#2^#3{\mathrel{\mathop{\kern 0pt#1}\limits_{#2}^{#3}}}

\begin{document}

\title[Phenomenologies of Hall-MHD turbulence]{Phenomenologies of intermittent Hall MHD turbulence}


\author [Mimi Dai]{Mimi Dai}

\address{Department of Mathematics, Stat. and Comp.Sci., University of Illinois Chicago, Chicago, IL 60607, USA}
\email{mdai@uic.edu}

\thanks{The author was partially supported by NSF grants DMS--1815069 and DMS--2009422.}





\begin{abstract}

We introduce the concept of intermittency dimension for the magnetohydrodynamics (MHD) to quantify the intermittency effect. With dependence on the intermittency dimension, we derive phenomenological laws for intermittent MHD turbulence with and without the Hall effect. In particular, scaling laws of dissipation wavenumber, energy spectra and structure functions are predicted. Moreover, we are able to provide estimates for energy spectra and structure functions which are consistent with the predicted scalings.   

\bigskip


\hspace{0.02cm}CLASSIFICATION CODE: 35Q35, 76D03, 76W05.
\end{abstract}

\maketitle

\section{Introduction}

The incompressible magnetohydrodynamics (MHD) with Hall effect featuring the physics of magnetic reconnection is governed by the system of partial differential equations (PDEs)
\begin{equation}\label{hmhd}
\begin{split}
\partial_tu+u\cdot\nabla u-B\cdot\nabla B+\nabla p=&\ \nu\Delta u,  \\
\partial_t B+u\cdot\nabla B-B\cdot\nabla u+d_i\nabla\times((\nabla\times B)\times B)=&\ \mu\Delta B, \\
\nabla\cdot u=&\ 0. 
\end{split}
\end{equation}
Here, $u$, $p$ and $B$ represent the fluid velocity field, scalar pressure, and magnetic field, respectively; they are unknown functions on the spacial-time domain $\Omega\times [0,\infty)$. 
The parameters $\nu$ and $\mu$ denote respectively the viscosity and resistivity.
The parameter $d_i$ stands for the ion inertial length, below the scale of which the ions tend to separate from the magnetic field.
Some simple facts follow from the general form (\ref{hmhd}):
\begin{itemize}
\item [(i)] If $B=0$, system (\ref{hmhd}) reduces to the Navier-Stokes equation (NSE). 
\item [(ii)] If $d_i=0$, system (\ref{hmhd}) reduces to the classical MHD, and the magnetic field is frozen into the fluid. 
\item [(iii)] If $d_i>0$, the Hall effect breaks the frozen-in property and the system can capture the fast magnetic reconnection process very well. 
\item [(iv)] At scales much smaller than $d_i$, the magnetic field is frozen again, albeit this time into the electron fluid.  In the limit of small scales, the ion flow appears too slow against the motion of electrons and tends to form a neutralizing background, i.e., $u$ vanishes. 
Thus (\ref{hmhd}) reduces to the electron magnetohydrodynamics (EMHD)
\begin{equation}\label{emhd}
\partial_tB+d_i\nabla\times((\nabla\times B)\times B)=\mu\Delta B.
\end{equation}
\item [(v)] The Hall term $\nabla\times ((\nabla\times B)\times B)$ appears to be more singular and exhibit higher order derivative than $u\cdot \nabla u$ and $B\cdot \nabla B$ in the NSE. It captures intricate dynamics responsible for striking turbulence phenomena and complicated energy cascade associated with magnetic reconnection processes. 
\item [(vi)] If $\nabla\cdot B(0)=0$ at the initial time, $\nabla\cdot B(t)=0$ remains for all the time $t>0$.
\end{itemize}
One observes immediately that the MHD system involves more complicated interactions than the NSE since the former contains nonlinear couplings of the magnetic field and the velocity field. The fact (v) indicates, from a surface level,  the dynamics of the Hall MHD is more intricate than that of the classical MHD.

\subsection{Background}

Proposed by Alfv\'en in 1942, the classical theory of MHD connects the Maxwell electrodynamics with the Navier-Stokes hydrodynamics. Over the next decades, classical MHD theory has evolved to lie at the heart of understanding most of the phenomena in plasma physics such as solar winds, interstellar clouds, planetary magnetospheres, etc. Derived from the classical MHD, Alfv\'en's famous frozen-in theorem demonstrates that magnetic field lines move with the ion flow. However, this frozen-in property is found to be invalid in some violent events, like solar flares. The widely accepted theory to explain the mystery is that the violent events involve active dynamics at small scales that are comparable or smaller than the ion inertial length $d_i$. At such small scales, ions tend to decouple from the magnetic field which becomes no longer frozen into the bulk plasma and changes topology through a rapid magnetic reconnection process. At the reconnection occurrence, an intense current sheet is created and a vast amount of energy gets released. Various models have been proposed to characterize the feature of the striking magnetic reconnection process. Among them, the MHD model with Hall effect is widely adapted, which is derived under the umbrella of the two-fluid reconnection theory. 

The MHD and Hall MHD models have been extensively studied by physicists \cite{Ber1, Ber2, BL1, Bha, Bis1, Bis2, Bold1, Bold2, BMC, CSM, Dav3, Dav4, DMV, Galtier, PPS, PFM, SchC} and mathematicians \cite{ADFL, CKS, CDL, CS, CWeng, CW, CDreg, Dai-uni, DuS, JO}. In particular,   the Hall MHD has attracted relentless interest in the community of mathematics in the past decades. Nevertheless, many peculiar phenomena in plasma physics remain to be resolved. Of crucial importance, the topic of turbulence is still an outstanding challenge in both of the  mathematics and physics communities in the new century. The (Hall) MHD turbulence plays a vital role in many complex plasma phenomena, such as the formation of accretion discs,  explosions on the surface of the Sun which lead to solar flares and coronal mass ejections, solar dynamo process, etc. These phenomena involve intricate interactions between the magnetic field and turbulent motions of the electrically conducting fluid, which cause complicated energy transfer between the kinetic and magnetic spectra. Applications of (Hall) MHD in  nuclear fusion and electrical power generation appear to be very inviting nowadays and require a thorough understanding of turbulence phenomena associated with the dynamics.


With the keenness of finding the order in chaos, numerous scientists have made extensive efforts to study the nature of turbulence in fluid motions for centuries. Our understanding of hydrodynamic turbulence has been greatly enriched since the middle of last century, with contributions of Kraichnan \cite{Kr1, Kr2}, Kolmogorov \cite{K41, K41b, K41c}, Onsager \cite{On}, 
Taylor \cite{Tay1, Tay2}, 
 etc. On one hand, the theory of hydrodynamic turbulence serves as a good foundation for building the MHD turbulence theory. On the other hand, in contrast with the NSE governing the hydrodynamics, the MHD models comprise richer nonlinear structures which are the origins of complicated energy cascade in turbulence phenomena. Especially, the Hall MHD model contains  a higher order nonlinear Hall term which is responsible for many striking turbulence phenomena. 
There is an enormous body of literature on theoretical, experimental and computer assisted studies of MHD turbulence. For MHD turbulence with low magnetic Reynolds number $R_m$, reasonable phenomenological models are available and well supported by numerical simulations. However, the theory of MHD turbulence with high $R_m$, which is a regime of great importance in plasma physics, is far from being satisfactory. In fact, there is still a lively debate over several existing phenomenological models for high-$R_m$ turbulence, which are not well consistent. Meanwhile, the understanding of Hall MHD turbulence with high $R_m$ is very limited, due to the difficulty that the PDE model (\ref{hmhd}) involves more intricate nonlinear interactions within a broader range of active space-time scales. 

Therefore, it urges advancements in theoretical study of the (Hall) MHD turbulence to formulate quantitative and testable predictions. In the present paper, ideas on quantifying intermittency effect will be emerged dimly into the study of phenomenological theory for both the Hall MHD and MHD systems. More details on the topic of intermittency will be provided in Subsection \ref{subsec-int} and Section \ref{sec-int}. 


\subsection{Review of K41 theory and K62 theory}
\label{sec-kol}
The celebrated Kolmogorov 1941 phenomenological turbulence theory \cite{K41} (referred as K41) for hydrodynamics was derived for homogeneous and isotropic flows under the assumption of self-similarity. It drew important predications on energy cascade from scale to scale, scaling laws for structure functions and energy spectrum, and even exact relation for the third order structure function. The main results comprise the four-fifths law for the third order structure function, two-thirds law for the second order structure function, negative five-thirds power law for the energy spectrum, and the derivation of Kolmogorov's dissipation wavenumber $\kappa_{\mathrm d}$ and dissipation scale $\eta$ with $\kappa_{\mathrm d}=1/\eta$, which separates the dissipation range from the inertial range. 
The following universality was postulated in the derivations: at very high Reynolds number, the small-scale statistical properties are uniquely and universally determined by the scale $\ell$, the mean energy dissipation rate $\varepsilon$ and the viscosity $\nu$. However, this universality was objected by Landau \cite{LL}. There is evidence that Kolmogorov was actually aware of such issue and modified the K41 theory in early 1960s \cite{K62}.
Nevertheless, as a consequence of K41 theory, investigations have been extended to the topics of estimating degrees of freedom, comparing macroscopic and microscopic length scales, finding the law of energy decay, etc.; and fruitful results have been established.

An important issue was that the plausible assumption of self-similarity (and homogeneity, isotropy) in K41 theory can be invalid for some turbulent flows. 
In fact, experimental evidences \cite{AGHA} show discrepancies from the K41 predictions in some situations. The deviation from K41 scaling suggests that small scales have fractal properties. As a general principal, the notion of intermittency was defined to correspond such deviation from the K41 theory. Two classes of phenomenological models with intermittency correction were introduced by Kolmogorov \cite{K62} to modify the K41 theory; the updated theory was referred as K62. In one class, intermittency is studied via velocity increments; while in another, intermittency is studied via dissipation fluctuation and a bridging argument of connecting inertial range quantities with such fluctuation, \cite{Fri, K62, Obu}. In the former class, based on the idea of quantifying space-filling of eddies in the Richardson cascade image, a fractal dimension parameter $D$ was introduced; and the so-called $\beta$-model, bifractal model, and multifractal model were derived to modify the K41 scaling laws with the dependence on $D$. In the latter class, the central idea was to define and measure multifractality in terms of the fluctuations of the local dissipation rather than in terms of velocity increments. In fact, the mean energy dissipation is a crucial quantity in the K41 theory; Landau's objection to the universality assumption also concerned mainly with the dissipation fluctuations. Kolmogorov \cite{K62} also found a bridge connecting the two multifractal formalizations. Predictions of the modified models in K62 theory by taking into account the intermittency effect are well consistent with experimental data, for instance, see \cite{MeS}.

\subsection{Review of scaling theories for MHD turbulence}
In the context of MHD turbulence, the magnetic field is a large-scale feature that remains at small scales \cite{NSch}. In contrast, large-scale features in hydrodynamic turbulence always regress to the Kolmogorov state at small scales. Therefore, MHD turbulence shares certain analogy with hydrodynamic turbulence but is primarily different from the latter. 
Systematic study of MHD turbulence was initiated by Iroshnikov \cite{Ir} and Kraichnan \cite{Kr1} in 1960s, who proposed a scaling theory (referred as IK) for flow with a uniform background magnetic field (Alfv\'en speed) and in a state of weak turbulence under the assumption of isotropy in the inertial range. By realizing the anisotropic feature of the turbulent motions along and across the background state, Goldreich and Sridhar \cite{GS1, GS2} proposed a phenomenological theory (noted as GS) standing on the critical balance conjecture in 1990s. As many physicists \cite{MGold, BL2} realized the crucial effect of alignment between the Els\"asser fields in the nonlinear interactions, Boldyrev \cite{Bold1, Bold2} presented an appealing theory based on the concept of dynamic alignment in 2000s. Brief review of the development of these theories will be highlighted below.

Let $B_0$ be the background uniform magnetic field. The Alfv\'en speed is the associated background velocity. 
Kraichnan \cite{Kr1} first realized that $B_0$ preserves at small scales through the energy cascade process. Relying on this argument, Iroshnikov and Kraichnan \cite{Ir, Kr1} derived the $-3/2$ power law for the magnetic energy spectrum. The scaling law was achieved via dimensional analysis and based on the assumption that turbulence in the inertial range is isotropic, in the spirit of the K41 theory. However, the fluctuations of the turbulent fields parallel and perpendicular to $B_0$ are not necessarily the same. In fact, since the magnetic field is frozen into the ion flow, the parallel variation goes with the propagation of Alfv\'en waves, and the perpendicular variation results from nonlinear interactions. It is thus natural to imagine that MHD turbulence is anisotropic at small scales; specifically, the characteristic scale of parallel fluctuation is larger than that of the perpendicular fluctuation. With such belief,  
Goldreich and Sridhar \cite{GS1, GS2} postulated that the parallel propagation wave period and the characteristic time of perpendicular variation are comparable, which is known as the critical balance conjecture. In view of the critical balance conjecture, Goldreich and Sridhar derived the $-5/3$ power law for the magnetic energy spectrum for the perpendicular fluctuation. This (perpendicular) scaling for anisotropic turbulence coincides with the K41 scaling for isotropic turbulence. While some solar wind turbulence observations were consistent with GS predictions, high resolution numerical simulations of MHD turbulence in early 2000s showed consistence with IK theory, see \cite{MGold, MBG}. 

By realizing the significance of the alignment between the two Els\"asser fields, Boldyrev made an assumption on the minimal degree of misalignment. Based on this alignment assumption, 
Boldyrev derived the $-3/2$ power law for the magnetic energy spectrum for the perpendicular variation. Thus, Boldyrev's perpendicular scaling returns to the IK thoery, which was supported by numerical studies \cite{MCB1, MCB2, PMBC1, PMBC2}.
Regarding the parallel cascade, it was derived \cite{GS2} that the magnetic energy spectrum in parallel direction satisfies the $-2$ power law, as a consequence of the critical balance conjecture. It remains true in the framework of Boldyrev's theory.

Nevertheless, up to the appearance of Boldyrev's theory, the story is not yet complete, as Beresnyak \cite{Ber1, Ber2, Ber3} brought up objection. In Beresnyak's opinion, 
Boldyrev's alignment theory violates the scaling symmetry of the Els\"asser fields equations and fails at small enough scales. 
Numerical simulations \cite{PMBC2, Ber3} performed by the two groups also showed disagreement.  Thereafter, serious debate on MHD turbulence follows from the competing theories and numerical results over the last decade. However, the intermittency effect brings hope to reconcile Boldyrev's alignment theory and Beresnyak's objection. The dynamic alignment theory can be interpreted as a qualitative ideology of intermittency. Recently,
it was realized in \cite{CSM, MSch} by the team of Schekochihin et al. that a model with intermittency consideration is crucial to derive the scaling laws. They revised Boldyrev's alignment theory by introducing the parallel outer scale as an extra parameter and incorporating dimensional correctness. The inclusion of the parallel outer scale indicates the invalidation of self-similarity of the MHD turbulence. They argue that the anisotropy depends on the local direction of the fluctuating fields. The turbulent field is then viewed as an ensemble of structures which have three scales corresponding to the parallel, perpendicular, and fluctuation directions. Analysis on the joint probability distribution of these quantities is carried through to fix scalings. Their revision of the alignment theory resembles the K62 theory for the hydrodynamic turbulence in which intermittency correction is the central idea to deal with the failure of self-similarity.

There is a vast literature of active research on other aspects of MHD turbulence. 
Without the intention of being  complete, we list a few topics. The magnetorotational dynamo was studied by statistical simulation \cite{SB}; 
the role of magnetorotational instability and plasmoid instability was investigated in \cite{CLHB}. The formation of accretion discs has been extensively studied by many scientists.
Another significant subject concerns ideal invariants for intermittent flows and their conservation versus anomalous dissipation \cite{BBV, Tay1, Tay2}. 

\subsection{Quantification of intermittency} 
\label{subsec-int}
Since the K62 theory, there has been growing interest in the study of intermittency. Various theoretical interpretations have been proposed, which are based on traditional statistical and probability theory of turbulence. 
In this paper, in order to quantify the intermittency effect, we will introduce an intermittency parameter -- intermittency dimension $\delta_b$ for the magnetic field, through the saturation level of Bernstein's inequality in harmonic analysis.  Bernstein's inequality provides quantitative relationships between different Lebegue norms associated with $L^p$ spaces. The essential idea is that the $L^p$ norms may be different for different values of $p$ for a turbulent field due to intermittency. According to our definition, when the $L^p$ norms of an $n$ dimensional field are the same for all $1\leq p\leq\infty$, the intermittency dimension is $n$; oppositely, when the $L^p$ norms are different up to the full saturation of Bernstein's inequality, the intermittency dimension is $0$. For instance, the Kolmogorov regime in K41 theory corresponds to the case of intermittency dimension being $3$, and the eddies at each scale occupy the whole region in the Richardson cascade image. The quantitative interpolation relationship equipped in Bernstein's inequality allows us to measure and quantify the intermittency effect of a turbulent field precisely by the intermittency dimension parameter. 

In this paper, we intend
to implement the quantification of intermittency effect into the MHD and Hall MHD turbulence theory. 
In particular, we will pursue the following objectives. 

\begin{itemize}
\item  Introduce magnetic intermittency dimension $\delta_{b}$ as a parameter to quantify the non-uniformity of a turbulent magnetic field; such concept of intermittency dimension will be adapted to turbulent Els\"asser fields as well;
\item   Establish phenomenological scaling laws of energy spectra and structure functions for the Hall MHD turbulence, with dependence on the intermittency dimension; and justify the phenomenology mathematically; 
\item  Explore transition scales between different energy cascade regimes with different energy spectra for the Hall MHD;
\item  Derive scaling laws of energy spectra and structure functions with intermittency correction for intermittent MHD turbulence; 
\item  Seek connections between findings above and the existing scaling theories of MHD turbulence, i.e. Iroshnikov-Kraichnan, Goldreich-Sridhar and Boldyrev's theories.
\end{itemize}


\section{Preliminaries and notations}
\label{sec-notations}

\subsection{Notations regarding constants}
Throughout the paper, we denote $A\lesssim B$ by an estimate of the form $A\leq cB$ for some constant $c$, and  $A\sim B$ an estimate of $c_1B\leq A\leq c_2B$ for constants $c_1$ and $c_2$.

\subsection{Littlewood-Paley decomposition}
The important parameter - intermittency dimension - will be defined via the saturation level of Bernstein's inequality. Therefore, we briefly introduce notations associated with Littlewood-Paley decomposition theory. Let $L$ be the domain size. We denote wavenumber $\lambda_q=\frac{2^q}{L}$ for integers $q$. For a tempered distribution vector field $v$ on $\mathbb T^3=[0,L]^3$, we denote $v_q$ by the $q$-th Littlewood-Paley projection of $v$. We also fix the notation
\[v_{\leq q}=\sum_{p\leq q} v_p, \ \ \ v_{> q}=\sum_{p> q} v_p.\]

\subsection{Energy flux}
For the EMHD (\ref{emhd}), we denote
\begin{equation}\notag
\Pi_{b,q}=d_i\int_{\mathbb T^3}((\nabla\times B)\times B) \cdot \nabla \times B_{<q}\, \mathrm{d} x
\end{equation}
by the magnetic energy flux below wavenumber $\lambda_q$. 
Thanks to the vector identity 
\[v\times w\cdot v=0,\]
we have
\begin{equation}\notag
\int_{\mathbb T^3}((\nabla\times B_{<q})\times B) \cdot \nabla \times B_{<q}\, \mathrm{d} x=0,
\end{equation}
which indicates that eddies larger than $\ell_q\sim \lambda_q^{-1}$ on average do not carry the energy across the scale $\ell_q$. Thus, in fact, we can write
\begin{equation}\notag
\Pi_{b,q}=d_i\int_{\mathbb T^3}((\nabla\times B_{\geq q})\times B) \cdot \nabla \times B_{<q}\, \mathrm{d} x.
\end{equation}
The formula can be further reduced to 
\begin{equation}\notag
\Pi_{b,q}=d_i\int_{\mathbb T^3}\sum_{p_1\geq q, p_2\geq q-1, p_3<q, |p_1-p_2|<2} ((\nabla\times B_{p_1})\times B_{p_2}) \cdot \nabla \times B_{p_3}\, \mathrm{d} x
\end{equation}
in view of the Fourier support of the projections. The idea is that remote scales do not contribute to the energy budget.
We also define the energy flux density as 
\begin{equation}\notag
\pi_{b,q}=d_i \sum_{p_1\geq q, p_2\geq q-1, p_3<q, |p_1-p_2|<2}((\nabla\times B_{p_1})\times B_{p_2}) \cdot \nabla \times B_{p_3}.
\end{equation}


\subsection{Energy spectrum and average energy dissipation rate}
We use the symbol $\left<\right>$ to represent time-space average of a vector field. 
For a vector field $v$, we denote 
\[\frac12\left<|\mathbb P_{\leq k}v|^2\right>\] by the mean energy per unit mass carried by wavenumber $\leq k$. The energy spectrum for $v$ can be defined as
\begin{equation}\notag
\mathcal E_{v}(k)=\frac12\frac{d}{dk} \left<|\mathbb P_{\leq k}v|^2\right>.
\end{equation}
Thus the total mean energy can be represented as 
\begin{equation}\notag
\frac12\left<|B|^2\right>=\int_0^\infty  \mathcal E_{v}(k)\,\mathrm{d}k.
\end{equation}
Note that we also have
\begin{equation}\notag
\mathcal E_{v}(\lambda_q)\sim \frac{\left<|v_q|^2\right>}{\lambda_q}.
\end{equation}
Denote the average dissipation rate of the magnetic energy by 
\[\varepsilon_{b}= \mu \left<\|\nabla B\|_{L^2}^2\right>.\]

\subsection{Structure function}
Let $\ell$ be the characteristic eddy size in the turbulence and $L$ the size of reference scale. Denote the characteristic scale of the fluctuating vector field $v$ by 
\[\delta_{\ell} v=v(x+\ell,t)-v(x,t).\]
The $p$-th order structure function of $v$ at scale $\ell$ is defined to be
 \[S_p(\ell)=\left<|\delta_{\ell} v|^p\right>.\]

\section{Intermittency dimension: a quantitative measure of intermittency effect}
\label{sec-int}

A central assumption of K41 turbulence theory is the self-similarity of the random velocity field in the inertial range. This assumption seems plausible but may well be invalid for intermittent flow. The revised K62 theory took into account the intermittency effect, which indicates certain fractal or multifractal properties of the turbulent flow. In K62 theory, intermittency was studied via velocity increments and dissipation fluctuation within the traditional framework of statistics and probability. Recently, intermittency was also analytically studied in \cite{CheS}. 

Different from the methods of K62 theory, the intermittency effect for the NSE/Euler flows is quantified in \cite{CDnse3-modes} by measuring the saturation level of Bernstein's inequality based on harmonic analysis techniques. We introduce the notion of intermittency dimension for a turbulent field as follows.

Let $L$ denote the length scale of the space domain $\Omega$. For integers $q\geq 0$, let $\lambda_q=2^q/L$ be the wavenumber of the $q$-th shell. For a vector field $u$, $u_q$ stands for the Littlewood-Paley projection of $u$ onto the $q$-th shell. One can understand it in the simple way that the Fourier transform of $u_q$ is supported on and near the $q$-th shell. 

\begin{Lemma}\label{le-bernstein}
[Bernstein's inequality]
Let $n$ be the spacial dimension and $p_2\geq p_1\geq1$. Then for all tempered distributions $u$, the inequality
\[\|u_q\|_{L^{p_2}}\leq C\lambda_q^{n(\frac1{p_1}-\frac1{p_2})}\|u_q\|_{L^{p_1}} \ \ \ \mbox{holds for a constant} \ \ C. \]
\end{Lemma}
A particular case of Bernstein's inequality for $n=3$, $p_2=\infty$ and $p_1=2$ reads as 
\[\lambda_q^{-1}\|u_q\|_{L^\infty}^2\leq C\lambda_q^2\|u_q\|_{L^2}^2.\]
On the other hand, we estimate the integral directly as
\[\int_{\Omega}|u_q|^2\, dx\leq |\Omega| \|u_q\|_{L^\infty}^2\]
which leads to 
\[c\lambda_0^3\|u_q\|_{L^2}^2\leq \|u_q\|_{L^\infty}^2\]
by realizing $\lambda_0=1/L$ and $ |\Omega|=\frac1{c} L^3$ for an absolute constant $c$. Combining the two inequalities above yields
\begin{equation}\label{bern1}
c\lambda_0^3\lambda_q^{-1}\|u_q\|_{L^2}^2\leq \lambda_q^{-1}\|u_q\|_{L^\infty}^2\leq C\lambda_q^2\|u_q\|_{L^2}^2.
\end{equation}
It is obvious that there is a scaling difference ($\lambda_q^3$) between the lower and upper bounds of the quantity $\lambda_q^{-1}\|u_q\|_{L^\infty}^2$.  Inspired by this observation, we introduce the intermittency parameter -- intermittency dimension, for a vector field $u$.

\begin{Definition}\label{int-def}
The intermittency dimension $\delta_u$ for a magnetic field $u(t)$ in three dimensions (3D) is defined as
\begin{equation}\label{def-int}
\delta_u:=\sup\left\{s\in\R:\left<\sum_{q}\lambda_q^{-1+s}\|u_q(t)\|_{L^\infty}^2\right>\leq C L^{-s}\left<\sum_{q}\lambda_q^2\|u_q(t)\|_{L^2}^2\right>\right\} 
\end{equation}
where the symbol $\left<\cdot\right>$ denotes time average and $C$ is a constant. 
\end{Definition} 
It follows from (\ref{bern1}) and (\ref{def-int}) that $\delta_{u}\in[0,3]$. Moreover, we have the following scaling relationship for a 3D vector field with intermittency dimension $\delta_u$
\begin{equation}\label{bern2}
\|u_q\|_{L^{p_2}}= C\lambda_q^{(3-\delta_u)(\frac1{p_1}-\frac1{p_2})}\|u_q\|_{L^{p_1}}
\end{equation}
for some constant $C$. We infer from (\ref{bern2}) that: 
\begin{itemize}
\item
when $\delta_u=3$, the $L^p$ norms for all $1\leq p\leq \infty$ are equivalent at each scale; 
\item
when $\delta_u=0$, the difference between different $L^p$ norms reaches the extreme scaling; in this case, Bernstein's inequality in Lemma \ref{le-bernstein} becomes an equality and hence is said to be saturated. 
\end{itemize}
According to this interpretation, the homogeneous, isotropic, and self-similar turbulent flow in K41 theory has intermittency dimension 3.  In the extreme case of intermittency dimension being 0, we understand it as there is only one eddy at each scale and the flow is very singular. On the other hand, numerical simulations and experimental studies show that $\delta\approx 2.7$ for classical hydrodynamics. By convention, a turbulent field with smaller intermittency dimension is said to be more intermittent. 

Intermittency dimension can be defined in the same way as in Definition \ref{int-def} for a magnetic field $B$ and the Els\"asser variables $Z^+=u+B$ and $Z^-=u-B$, respectively, denoted by $\delta_b$, $\delta^+$ and $\delta^-$. Analogously, we have the scaling relationship 
\begin{equation}\label{bern-intermit}
\begin{split}
\|B_q\|_{L^{p_2}}=&\ C\lambda_q^{(3-\delta_b)(\frac1{p_1}-\frac1{p_2})}\|B_q\|_{L^{p_1}},\\
\|Z^+_q\|_{L^{p_2}}=&\ C\lambda_q^{(3-\delta^+)(\frac1{p_1}-\frac1{p_2})}\|Z^+_q\|_{L^{p_1}},\\
\|Z^-_q\|_{L^{p_2}}=&\ C\lambda_q^{(3-\delta^-)(\frac1{p_1}-\frac1{p_2})}\|Z^-_q\|_{L^{p_1}}.
\end{split}
\end{equation}

A statistical concept of intermittency is also introduced in the context of stochastic processes by Khoshnevisan \cite{Kho}. The definition of \cite{Kho} is not quantitative but rather qualitative. That is, it can tell whether a random field is intermittent or not; it can not describe how intermittent the random field is. Nevertheless, definition (\ref{def-int}) is quantitative; it measures how intermittent a turbulent field is.

\section{Phenomenologies of intermittent turbulence for Hall MHD and estimates} 
\label{sec-hmhd}

In this section, we aim to develop some phenomenological scaling laws of turbulence for the MHD with Hall effect by including intermittency dimensions of the velocity and the magnetic field. The Hall MHD turbulence is is not well understood, by virtue of the intricate coupling of fluid velocity and magnetic field, and the extra complexity brought in by the Hall effect. In particular, the Hall term launches new physics into the system at small scales, which cause the system more ``chaotic". We will emphasize on predicting scaling laws for energy spectrum and structure functions; we will also extend the study to find the transition scales separating different regimes of energy cascade.

\subsection{Phenomenologies of intermittent EMHD and estimates}
For the EMHD model (\ref{emhd}) with reduced complications of multi-scales and nonlinear interactions, we have the main scaling laws regarding the dissipation wavenumber that separates the dissipation range from the sub-ion range, magnetic energy spectrum in the sub-ion range, and structure functions.


\begin{Conjecture}\label{thm-dissipation-emhd}
{\it Let $B$ be a solution of the 3D EMHD (\ref{emhd}) with intermittency dimension $\delta_b$. 
There exists a dissipation wavenumber $\kappa_{\mathrm d}^{e}$ for the 3D EMHD (\ref{emhd}) with the scaling }
\begin{equation}\label{hkd}
\kappa_{\mathrm d}^e\sim \left(\mu^{-3}d_i^2\varepsilon_{ b}\right)^{\frac1{\delta_{b}-1}},
\end{equation}
{\it such that $\kappa_{\mathrm d}^{e}$ separates the dissipation range from the sub-ion range. }
\end{Conjecture}

\begin{Conjecture}\label{q41}
{\it For the 3D EMHD equation (\ref{emhd}), the magnetic energy spectrum $\mathcal E_{b}(k)$ in the sub-ion range obeys the scaling law }
\begin{equation}\label{hen-in}
\mathcal E_{b}(k)\sim \left(d_i^{-1}\varepsilon_{b}\right)^{\frac23}k^{\frac{\delta_{b}-10}3}.
\end{equation}
\end{Conjecture}
\begin{Remark}
If the magnetic field is homogeneous, isotropic, and self-similar, i.e. in the case of $\delta_b=3$, it follows from (\ref{hen-in}) that 
\[\mathcal E_{b}(k)\sim \left(d_i^{-1}\varepsilon_{b}\right)^{\frac23}k^{\frac73}\]
which coincides with the scaling derived by physicists, for instance, see \cite{Galtier}.
\end{Remark}

We point out that Conjecture \ref{thm-dissipation-emhd} is in analogy with Kolmogorov's prediction on the dissipation wavenumber for hydrodynamics with intermittency correction. In fact, Kolmogorov's prediction on the dissipation wavenumber was recently justified mathematically by the author and collaborator in \cite{CDnse3-modes, CDK}. We expect to be able to justify Conjecture \ref{thm-dissipation-emhd} as well following the line of the previous work \cite{CDnse3-modes}, and we will address it in future work.

A rigorous proof of Conjecture \ref{q41} is unlikely to be achieved with existing techniques. Nevertheless, we are able to obtain upper and lower bounds for the magnetic energy spectrum. Namely, we have:

\begin{Theorem}\label{thm-bounds-spectrum}
Denote 
\[\bar\varepsilon_b=\sup_{q}d_i\lambda_q^2\left<|B_q|^3\right>, \ \ \ \ \underline\varepsilon_b=\inf_{q} \left<|\pi_{b,q}|\right>.\]
The energy spectrum $\mathcal E_{b}(k)$ satisfies the following upper bound
\[\mathcal E_{b}(k)\lesssim (d_i^{-1}\bar\varepsilon) ^{\frac23} k^{-\frac73} (Lk)^{\frac{\delta_b}{3}-1}\]
and an average lower bound
\begin{equation}\notag
\sum_{p} K_{q-p}^{\frac23} \lambda_p^{\frac13(10-\delta_b)} \mathcal E_b(\lambda_p)\geq (d_i^{-1}\underline\varepsilon) ^{\frac23}
\end{equation}
for any $q$, with $K_q=\lambda_{|q|}^{-\frac13}$.
\end{Theorem}

Beside dissipation wavenumber and energy spectrum, we can also predict scaling law for structure functions of the intermittent magnetic field. 
Indeed, scaling analysis suggests:





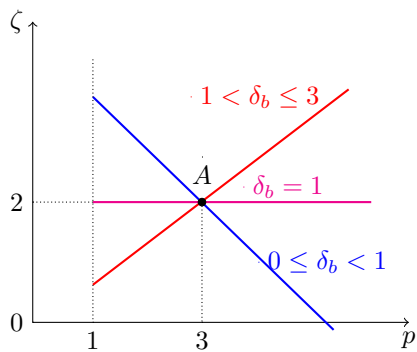
\begin{wrapfigure}[14]{r}{0.4\textwidth}
\begin{tikzpicture}
\draw [<->] (0,4) node[left]{$\zeta$}  --(0,0)  -- (5,0) node[below]{$p$};
\draw  (0.8, 0) -- node[below]{$1$} (0.8, 0); 
\draw  (0, 0) -- node[left]{$0$} (0, 0); 
\draw  (0, 1.6) -- node[left]{$2$} (0, 1.6);
\draw  [thick, magenta] (0.8,1.6)  -- (4.5,1.6) ; 
\draw  [thick, blue] (0.8,3)  -- (4,-0.1) ;  
\draw  [thick, red] (0.8,0.5)  -- (4.2,3.1) ;  
\draw [densely dotted] (0.8,3.5)  -- (0.8,0) ;  
\draw [densely dotted] (0.8,1.6)  -- (0,1.6) ; 
\draw [densely dotted] (2.25,1.6)  -- (2.25,0) ; 
\draw  (2.25, 0) -- node[below]{$3$} (2.25, 0); 
\draw  (2.25, 2.2) -- node[below]{$A$} (2.25, 2.2); 
\draw [red] (2.1,3)node[right]{$1<\delta_b\leq 3$} -- (2.1,3)  ;
\draw  [magenta] (2.8,1.8)node[right]{$\delta_b=1$} -- (2.8,1.8)  ;
\draw [blue] (3,0.8)node[right]{$0\leq\delta_b< 1$} -- (3,0.8)  ;
\draw [fill] (2.25,1.6) circle [radius=0.05];
\end{tikzpicture}
\caption{Structure function exponent as a function of $p$ for EMHD with different intermittency level.}\label{fig3}
\end{wrapfigure}

\begin{Conjecture}\label{q42}
{\it
Let $B$ be a solution of the 3D EMHD (\ref{emhd}) with intermittency dimension $\delta_b$. The $p$-th order structure function  $S_p(\ell)=\left<|\delta_{\ell} B|^p\right>$ has the scaling  }
\begin{equation}\label{sf2}
S_p(\ell)\sim \left(d_i^{-1}\varepsilon_{b}\right)^{\frac p3}\ell^{\frac{2p}{3}+(3-\delta_b)(1-\frac{p}{3})}.
\end{equation}
\end{Conjecture}

We denote the exponent of the structure function scaling by 
$\zeta(p,\delta_b)=\frac{2p}{3}+(3-\delta_b)(1-\frac{p}{3})$.

\begin{Remark}\label{rk-structure1}
We notice that $\zeta(3,\delta_b)=2$ which does not dependent on the intermittency dimension. It suggests that there is an exact law for the 3rd-order structure function.
\end{Remark}
\begin{Remark}\label{rk-structure2}
We also notice that $\zeta(p,\delta_b)$ is a linear function in both $p$ and $\delta_b$.
\end{Remark}

For homogeneous, isotropic, and self-similar magnetic field, that is $\delta_b=3$,  the scaling (\ref{sf2}) represents $4/3$ law for the second order structure function, which is consistent with the derivation of physicists, for instance, see \cite{Galtier}.  Two special cases are: if $p=3$, the scaling of the third structure function is of $\ell^{2}$ for the magnetic field with any intermittency dimension, see point A in Figure \ref{fig3}; if $\delta_b=1$, the scaling is of $\ell^{2}$ for the structure function of any order. It is also important to notice that how intermittency level affects the property of structure functions. The exponent $\zeta(p,\delta_b)$ increases with $p$ if  $\delta_b>1$ (red line in Figure \ref{fig3}); while it decreases with $p$ if  $\delta_b<1$ (blue line in Figure \ref{fig3}). The second and third order structure functions are illustrated for homogeneous isotropic self-similar turbulence and extremely anisotropic turbulence respectively in Figure \ref{fig1} and Figure \ref{fig2}.

\begin{minipage}{.48\textwidth}
\begin{center}
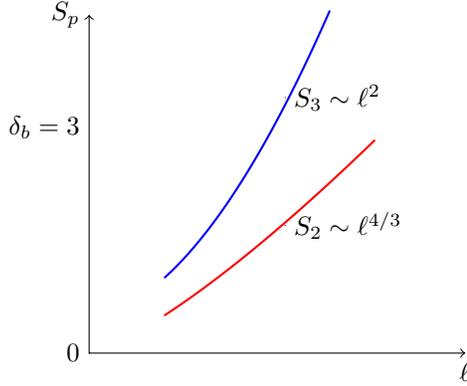

\begin{tikzpicture}
\draw [<->] (0,4.5) node[left]{$S_p$}  --(0,0)  -- (5,0) node[below]{$\ell$};
\draw  (0,3) node[left]{$\delta_b=3$}; 
\draw  (0, 0) -- node[left]{$0$} (0, 0); 
\draw [thick, red] plot[smooth,domain=1:3.8] (\x, {(0.5)*\x^1.3});
\draw [thick, blue] plot[smooth,domain=1:3.2] (\x, {(0.5)*\x^1.8+0.5});    
\draw  (2.6,3.4)node[right]{$S_3\sim \ell^2$} -- (2.6,3.4)  ;
\draw  (2.6, 1.7)node[right]{$S_2\sim \ell^{4/3}$} -- (2.6,1.7)  ;
\end{tikzpicture}
\captionof{figure}{Second (red) and third (blue) order structure functions for homogeneous isotropic self-similar EMHD turbulence.}
\label{fig1}
\end{center}
\end{minipage}
\begin{minipage}{.48\textwidth}
\begin{center}
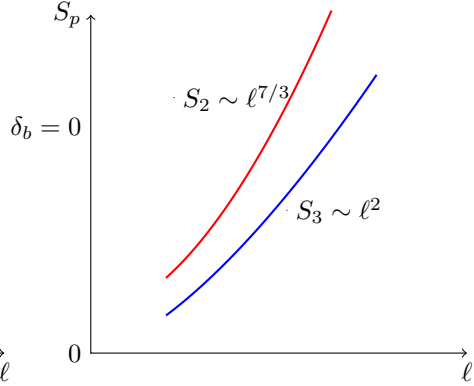

\begin{tikzpicture}
\draw [<->] (0,4.5) node[left]{$S_p$}  --(0,0)  -- (5,0) node[below]{$\ell$};
\draw  (0,3) node[left]{$\delta_b=0$}; 
\draw  (0, 0) -- node[left]{$0$} (0, 0); 
\draw  [thick, blue] plot[smooth,domain=1:3.8] (\x, {(0.5)*\x^1.5});
\draw  [thick, red] plot[smooth,domain=1:3.2] (\x, {(0.5)*\x^1.8+0.5});    
\draw  (1.1,3.4)node[right]{$S_2\sim\ell^{7/3}$} -- (1.1,3.4)  ;
\draw  (2.6, 1.9)node[right]{$S_3\sim\ell^2$} -- (2.6,1.9)  ;
\end{tikzpicture}
\captionof{figure}{Second (red) and third (blue) order structure functions for extremely anisotropic EMHD turbulence.}
\label{fig2}
\end{center}
\end{minipage}


Although lack of a proof of Conjecture \ref{q42}, an upper bound on the structure functions with $2\leq p\leq3$ can be established.

\begin{Theorem}\label{thm-sf-bound}
Assume $\delta_b\in[1,3]$. Let $2\leq p\leq 3$. There exists a constant $C_p>0$ such that  
\[S_p(\ell)\leq C_p \left(d_i^{-1}\bar\varepsilon_{b}\right)^{\frac p3}\ell^{\frac{2p}{3}+(3-\delta_b)(1-\frac{p}{3})}.\]
\end{Theorem}

\subsection{Phenomenologies of intermittent Hall MHD}

The situation for the Hall MHD system is more complicated than that of the EMHD, due to more intricate nonlinear couplings and interactions. Moreover, the intermittency dimension $\delta_u$ of the velocity field and $\delta_b$ of the magnetic field both play vital roles here.  
The interesting question bas been raised in the field: whether the fluid velocity $u$ or the magnetic field $B$ plays a dominant role in the dynamics? We expect that the answer depends on the intermittency level of both $\delta_u$ and $\delta_b$. The scaling laws of intermittent Hall MHD turbulence should also depend on $\delta_u$ and $\delta_b$.

Denote the average dissipation rate of the kinetic
energy by $\varepsilon_{u}= \nu \left<\|\nabla u\|_{L^2}^2\right>$ and the kinetic energy 
spectrum by $\mathcal E_u(k)$. We consider the situation that 
 the velocity plays a dominant role, i.e., the influence of the velocity field over the magnetic field is stronger than the influence of the magnetic field over the velocity field. In this case, it is natural to assume that the velocity field is more intermittent than the magnetic field, i.e. $\delta_u\leq \delta_b$. We have the following predictions in this regime.
\begin{Conjecture}\label{q43}
{\it
Assume $\delta_u\leq \delta_b$. In the kinetic inertial range, the kinetic energy spectrum of the intermittent Hall MHD system exhibits the scaling 
\[\mathcal E_u(k)\sim \varepsilon_u^{\frac23}k^{\frac{\delta_u-8}{3}}.\]
The magnetic energy spectrum exhibits two power laws,
\begin{equation}\notag
\mathcal E_b(k)\sim
\begin{cases}
\varepsilon_b^{\frac12}\varepsilon_u^{\frac16}k^{\frac{\delta_u+3\delta_b}{12}-\frac83}, \ \ \ \mbox{ in ion-inertial range,}\\
(d_i^{-1}\varepsilon_b)^{\frac23} k^{\frac{\delta_b-10}{3}}, \ \ \mbox{ in sub-ion range}.
\end{cases}
\end{equation}
}
\end{Conjecture}

We proceed to further predict the transition scales that separate different regimes of cascade.

\begin{Conjecture}\label{q44}
{\it
The kinetic dissipation wavenumber $\kappa_{\mathrm d}^u$ that separates the dissipation range from the kinetic inertial range for the fluid velocity has the scaling
\begin{equation}\label{hmhd-kd-u}
\kappa_{\mathrm d}^u\sim \left(\nu^{-3}\varepsilon_{u}\right)^{\frac1{\delta_{u}+1}}.
\end{equation}
For the magnetic field, the critical wavenumber $\kappa_{\mathrm {i}}^b$ that separates the ion-inertial range from the sub-ion range satisfies 
\[\kappa_{\mathrm {i}}^b\sim d_i^{-1};\]
and the magnetic dissipation wavenumber $\kappa_{\mathrm d}^b$ that separates the sub-ion range from the dissipation range obeys the scaling
\begin{equation}\label{hmhd-kd-b}
\kappa_{\mathrm d}^b\sim \left(\mu^{-3}\varepsilon_{u}^{-\frac12}\varepsilon_{ b}^{\frac32}\right)^{\frac{1}{\frac{3\delta_b+\delta_u}{4}+1}}.
\end{equation}
 }
\end{Conjecture}

From (\ref{hmhd-kd-u}) and (\ref{hmhd-kd-b}), one can observe that larger intermittency dimensions $\delta_u$ and $\delta_b$  indicate higher regularity and henceforth smaller dissipation wavenumber $\kappa_{\mathrm d}^u$ and $\kappa_{\mathrm d}^b$ and narrower kinetic inertial range and magnetic sub-ion range.  As a principal application, the dissipation wavenumber is often used to estimate the number of degrees of freedom and thereby brings hope to improve numerical simulation algorithms.

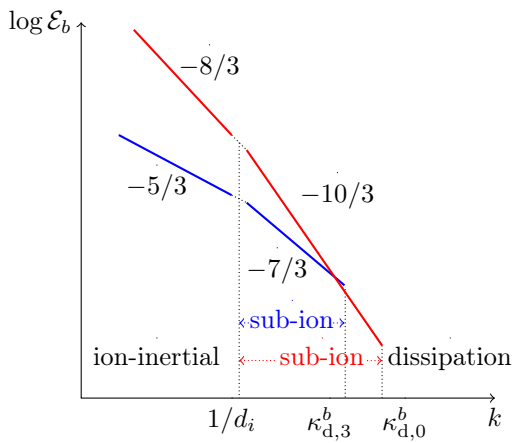
\begin{wrapfigure}[14]{r}{0.5\textwidth}
\begin{tikzpicture}
\draw [<->] (0,5) node[left]{$\log\mathcal E_b$}  --(0,0)  -- (5.5,0) node[below]{$k$};
\draw  [thick, blue] (0.5,3.5)  -- (2,2.7) ; 
\draw  [thick, blue] (2.2,2.6)  -- (3.5,1.5) ; 
\draw [densely dotted] (2, 2.7)  -- (2.2,2.6) ; 
\draw [densely dotted] (2.1, 3.4)  -- (2.1,0) ; \draw [densely dotted] (3.51, 1.45)  -- (3.51,0) ; 
\draw [thick, red] (0.7,4.9)  -- (2, 3.5) ;  \draw [thick, red] (2.2, 3.3)  -- (4, 0.7) ;  \draw [densely dotted] (2, 3.5)  -- (2.2,3.3) ; 
\draw [densely dotted] (4, 0.7)  -- (4,0) ;  
\draw  (2,0) -- node[below]{$1/d_i$}  (2,0)  ;
\draw  (3.3,0) -- node[below]{$\kappa_{\mathrm d,3}^b$}  (3.3,0)  ;
\draw  (4.3,0) -- node[below]{$\kappa_{\mathrm d,0}^b$}  (4.3,0)  ;
\draw  (1,0.8) -- node[below]{ion-inertial}  (1,0.8)  ;
\draw [blue] (2.8,1.3) -- node[below]{sub-ion}  (2.8,1.3)  ;
\draw [blue, densely dotted] [<-](2.1, 1) -- (2.25, 1);
\draw [blue, densely dotted] [->](3.3, 1) -- (3.5, 1);
\draw [red] (3.2,0.8) -- node[below]{sub-ion}  (3.2,0.8)  ;
\draw [red, densely dotted] [<-](2.1, 0.5) -- (2.6, 0.5);
\draw [red, densely dotted] [->](3.8, 0.5) -- (4, 0.5);
\draw  (4.9,0.8) -- node[below]{dissipation}  (4.9,0.8)  ;
\draw  (1,3.1) -- node[below]{$-5/3$}  (1,3.1)  ;
\draw  (2.6,2) -- node[below]{$-7/3$}  (2.6,2)  ;
\draw  (1.7,4.7) -- node[below]{$-8/3$}  (1.7,4.7)  ;
\draw  (3.4,3) -- node[below]{$-10/3$}  (3.4,3)  ;
\end{tikzpicture}
\caption{Magnetic energy spectra of Hall MHD when $\delta_u=\delta_b=3$ (blue lines) and when $\delta_u=\delta_b=0$ (red lines).}
\label{fig4}
\end{wrapfigure}

Figure \ref{fig4} illustrates the information contained in Conjecture \ref{q43} and Conjecture \ref{q44}.
We notice that the kinetic energy spectrum does not depend on the magnetic field intermittency dimension, which is consistent with the hypothesis that the velocity plays a dominant role. One can also see that the magnetic energy spectrum in ion-inertial range depends on average dissipation rates of both the kinetic energy and magnetic energy, and intermittency effect of both the velocity and magnetic field. It indicates intricate interactions and coupling within the ion-\\
inertial range. While in the sub-ion \\
range, the magnetic energy spectrum \\
does not depend on the velocity field, \\
since the Hall term plays a dominant \\
role within this regime.
The special case of $\delta_u=\delta_b=3$ for homogeneous isotropic self-similar turbulence has the scaling 
$k^{-7/3}$
 in the ion-inertial range and $k^{-5/3}$ in sub-
 ion range, which again is consistent with physics 
 phenomenology \cite{Galtier}. While the most steep scaling is that of the 
most extremely anisotropic turbulence with $\delta_u=\delta_b=0$.

Similarly as for Conjecture \ref{thm-dissipation-emhd}, it is hopeful to justify Conjecture \ref{q44} by applying the framework of determining wavenumber and the wavenumber splitting approach developed in the previous work \cite{CDnse3-modes} for 3D Navier-Stokes equation. While a rigorous proof will be pursued in future research, a heuristic analysis will be provided in Section \ref{sec-justification-hmhd} to motivate the scalings.

By the definition of structure function, we can see that the structure function of a vector field does not depend on the particular equation it satisfies, rather it depends on the intermittency effect. Thus, in view of Conjecture \ref{q42} and Theorem \ref{thm-sf-bound}, we state the scaling law of structure functions and their estimate for the Hall MHD system (\ref{hmhd}) in the following.

\begin{Conjecture}\label{q45}
{\it
Let $(u, B)$ be a solution of the 3D Hall MHD (\ref{hmhd}) with intermittency dimension $(\delta_u, \delta_b)$. The $p$-th order structure functions  $S_{b,p}=\left<|\delta_{\ell} B|^p\right>$ and $S_{u,p}=\left<|\delta_{\ell} u|^p\right>$ satisfy  }
\begin{equation}\label{sf-hmhd}
\begin{split}
S_{b,p}(\ell)\sim &\ \left(d_i^{-1}\varepsilon_{b}\right)^{\frac p3}\ell^{\frac{2p}{3}+(3-\delta_b)(1-\frac{p}{3})},\\
S_{u,p}(\ell)\sim &\ \left(\varepsilon_{u}\right)^{\frac p3}\ell^{\frac{p}{3}+(3-\delta_u)(1-\frac{p}{3})}.
\end{split}
\end{equation}
\end{Conjecture}

\begin{Theorem}\label{thm-sf-bound-hmhd}
Let $(u, B)$ be a solution of the 3D Hall MHD (\ref{hmhd}) with intermittency dimension $(\delta_u, \delta_b)$. 
Assume $\delta_u\in[0,3]$ and $\delta_b\in[1,3]$. Let $2\leq p\leq 3$. There exists a constant $C_{p}>0$ such that  
\begin{equation}\notag
\begin{split}
S_{b,p}(\ell)\leq &\ C_p \left(d_i^{-1}\bar\varepsilon_{b}\right)^{\frac p3}\ell^{\frac{2p}{3}+(3-\delta_b)(1-\frac{p}{3})},\\
S_{u,p}(\ell)\leq &\ C_p\left(\varepsilon_{u}\right)^{\frac p3}\ell^{\frac{p}{3}+(3-\delta_u)(1-\frac{p}{3})}.
\end{split}
\end{equation}
\end{Theorem}

\begin{Remark}
In the case $\delta_b< \delta_u$, similar phenomenologies as in Conjectures \ref{q43}-\ref{q45} can be derived as well. In this situation, one has to analyze which nonlinear interactions dominate in the inertial range and which nonlinear terms in (\ref{scal-u1})-(\ref{scal2}) to be used at each step of the scaling analysis.
\end{Remark}


\section{Derivation of phenomenologies and proof of estimates}
\label{sec-justification-hmhd}

\subsection{Scaling derivation of Conjecture \ref{thm-dissipation-emhd} and Conjecture \ref{q41}}
\label{sec-proof-con41}
The heuristic scaling analysis starts with the formal energy identity for the EMHD equation (\ref{emhd}),
\begin{equation}\notag
\frac12\frac{d}{dt}\int_{\mathbb T^3}|B|^2\, \mathrm{d}x+d_i\int_{\mathbb T^3}\nabla\times((\nabla\times B)\times B)\cdot B\, \mathrm{d}x+\mu \int_{\mathbb T^3}|\nabla B|^2\, \mathrm{d}x=0.
\end{equation}
The energy identity suggests that the nonlinear flux $d_i\|\nabla\times((\nabla\times B)\times B)\cdot B\|_{L^1}$ and the dissipation term $\mu\|\nabla B\|_{L^2}^2$ have the same scaling as the dissipation rate of the magnetic energy $E_b(t)=\frac12\|B(t)\|_{L^2}^2$. We consider such scaling relationship at the level of the $q$-th shell. In fact, projecting the equation (\ref{emhd}) onto the $q$-th shell, taking inner product with $B_q$ and integrating over $\mathbb T^3$ gives us the energy law
\begin{equation}\label{eq-energy-emhd}
\frac12\frac{d}{dt}\int_{\mathbb T^3}|B_q|^2\, \mathrm{d}x+d_i\int_{\mathbb T^3}\nabla\times((\nabla\times B)\times B)_q\cdot B_q\, \mathrm{d}x+\mu \int_{\mathbb T^3}|\nabla B_q|^2\, \mathrm{d}x=0.
\end{equation}

Denote $\varepsilon_{b,q}=\mu\|\nabla B_q\|_{L^2}^2$, which can represent the energy dissipation rate at the $q$-th shell, since it follows from (\ref{eq-energy-emhd}) that
\[\frac12\frac{d}{dt}\int_{\mathbb T^3}|B_q|^2\, \mathrm{d}x\sim \mu\|\nabla B_q\|_{L^2}^2.\]
Further, we can infer the following scaling relationship from (\ref{eq-energy-emhd}), by applying H\"older's inequality and the saturated Bernstein's relationship (\ref{bern-intermit})
\begin{equation}\notag
\begin{split}
\mu \lambda_q^2\|B_q\|_{L^2}^2\sim \mu\|\nabla B_q\|_{L^2}^2\sim &\ d_i\|\nabla\times((\nabla\times B_q)\times B_q)\cdot B_q\|_{L^1}\\
\sim &\ d_i\lambda_q^2\|B_q\|_{L^2}^2\|B_q\|_{L^\infty}\\
\sim &\ d_i\lambda_q^{\frac{7-\delta_b}{2}}\|B_q\|_{L^2}^3.
\end{split}
\end{equation}
It follows that 
\begin{equation}\label{scaling-bql2}
\|B_q\|_{L^2}\sim \mu d_i^{-1} \lambda_q^{\frac{\delta_b-3}{2}}.
\end{equation}
Applying (\ref{scaling-bql2}) to the energy dissipation rate of the $q$-th shell, we obtain
\begin{equation}\label{scaling-vareq}
\varepsilon_{b,q}=\mu\|\nabla B_q\|_{L^2}^2 \sim \mu \lambda_q^2\|B_q\|_{L^2}^2\sim \mu^3 \lambda_q^{\delta_b-1}d_i^{-2}.
\end{equation}
We extract the scaling relationship from (\ref{scaling-vareq})
\begin{equation}\notag
\lambda_q\sim (\mu^{-3}d_i^2\varepsilon_{b,q})^{\frac{1}{\delta_b-1}}
\end{equation}
which motivates the scaling law (\ref{hkd}) in Conjecture \ref{thm-dissipation-emhd}.

On the other hand, (\ref{scaling-vareq}) also implies 
\begin{equation}\label{est-mu-varep}
\mu\sim d_i^{\frac23}\varepsilon_{b,q}^{\frac13} \lambda_q^{\frac13(1-\delta_b)}.
\end{equation}
Therefore, combining (\ref{scaling-bql2}) and (\ref{est-mu-varep}), we infer
\begin{equation}\notag
\mathcal E_b(\lambda_q)\sim \frac{\left<|B_q|^2\right>}{\lambda_q}\sim \mu^2d_i^{-2}\lambda_q^{\delta_b-4}\sim (d_i^{-1}\varepsilon_{b,q})^{\frac23}\lambda_q^{\frac{\delta_b-10}{3}}
\end{equation}
which suggests the energy spectrum scaling (\ref{hen-in}) of Conjecture \ref{q41}.

\subsection{Proof of Theorem \ref{thm-bounds-spectrum}.}
Since $\mathcal E_b(\lambda_q)\sim \frac{\left<|B_q|^2\right>}{\lambda_q}$, let us first estimate $\left<|B_q|^2\right>$. By the definition of structure function, we have, applying H\"older's inequality
\begin{equation}\notag
\begin{split}
\left<|B_q|^2\right>=&\ \frac{1}{T|\Omega|}\int_0^T\int_\Omega |B_q|^2\,\mathrm{d}x\mathrm{d}t\\
\leq &\ \frac{1}{T|\Omega|}\int_0^T\left(\int_\Omega |B_q|^3\,\mathrm{d}x\right)^{\frac23}\mathrm{d}t \cdot V_q^{\frac13}\\
\leq &\ \frac{1}{T|\Omega|}\left(\int_0^T\int_\Omega |B_q|^3\,\mathrm{d}x\mathrm{d}t\right)^{\frac23} \cdot T^{\frac13}V_q^{\frac13}\\
\lesssim &\ \left<|B_q|^3\right>^{\frac23} V_q^{\frac13}
\end{split}
\end{equation}
where $V_q=(\ell_q/L)^{3-\delta_b}$ is the active volume at scale $\ell_q\sim\lambda_q^{-1}$. Recall
\[\bar\varepsilon_b=\sup_{q}d_i\lambda_q^2\left<|B_q|^3\right>.\]
Thus, we have
\begin{equation}\label{est-bq-l2-average}
\left<|B_q|^2\right>\leq \left<|B_q|^3\right>^{\frac23} (\ell_q/L)^{\frac{3-\delta_b}{3}}\leq (d_i^{-1}\bar\varepsilon)^{\frac23} \lambda_q^{-\frac43}(\ell_q/L)^{\frac{3-\delta_b}{3}}.
\end{equation}
It follows from (\ref{est-bq-l2-average}) that
\begin{equation}\notag
\frac{\left<|B_q|^2\right>}{\lambda_q}\leq (d_i^{-1}\bar\varepsilon)^{\frac23} \lambda_q^{-\frac73}(\ell_q/L)^{\frac{3-\delta_b}{3}}\lesssim (d_i^{-1}\bar\varepsilon)^{\frac23} \lambda_q^{-\frac73}(\lambda_qL)^{\frac{\delta_b}{3}-1}
\end{equation}
which gives the upper bound of the energy spectrum.

We move forward to establish the lower bound of the energy spectrum.
By (\ref{bern-intermit}), we have
\[\|B_q\|_{L^3}\sim \lambda_q^{\frac16(3-\delta_b)}\|B_q\|_{L^2},\]
and hence
\begin{equation}\label{est-bq-l2-l3}
\left<|B_q|^3\right>\sim  \lambda_q^{\frac12(3-\delta_b)} \left<|B_q|^2\right>^{\frac32}.
\end{equation}
Applying H\"older's inequality, we obtain
\begin{equation}\notag
\begin{split}
\left|\int_{\mathbb T^3}\pi_{b,q}\,\mathrm{d}x\right|\leq&\ d_i \sum_{p_1\geq q-1, p_2\geq q, p_3<q, |p_1-p_2|<2} \int_{\mathbb T^3}\left| ((\nabla\times B_{p_1})\times B_{p_2})\cdot \nabla\times B_{p_3}\right|\, \mathrm{d}x\\
\leq&\ d_i \sum_{p_1\geq q-1, p_2\geq q, p_3<q, |p_1-p_2|<2}\lambda_{p_1}\|B_{p_1}\|_{L^3} \|B_{p_2}\|_{L^3} \lambda_{p_3}\|B_{p_3}\|_{L^3}\\
\lesssim &\ d_i \sum_{p_1> q-2, p_3<q}\lambda_{p_1}\|B_{p_1}\|_{L^3}^2 \lambda_{p_3}\|B_{p_3}\|_{L^3}\\
\leq&\ d_i \left( \sum_{p> q-2}\lambda_{p}^{\frac32}\|B_{p}\|_{L^3}^3 \right)^{\frac23}
\left( \sum_{p< q}\lambda_{p}^{3}\|B_{p}\|_{L^3}^3 \right)^{\frac13}.
\end{split}
\end{equation}
Rearranging the wavenumber multiples in the last inequality, we have
\begin{equation}\label{est-pi-bq}
\begin{split}
\left|\int_{\mathbb T^3}\pi_{b,q}\,\mathrm{d}x\right|
\lesssim &\ d_i \left( \sum_{p> q-2}\lambda_{p}^{-\frac12}\lambda_q^{\frac12}\lambda_p^2\|B_{p}\|_{L^3}^3 \right)^{\frac23}
\left( \sum_{p< q}\lambda_{p}\lambda_{q}^{-1}\lambda_{p}^{2}\|B_{p}\|_{L^3}^3 \right)^{\frac13}\\
\lesssim &\ d_i \sum_{p} K_{q-p} \lambda_{p}^{2}\|B_{p}\|_{L^3}^3
\end{split}
\end{equation}
with the kernel $K_q=\lambda_{|q|}^{-\frac13}$.
Recall that \[\underline\varepsilon_b=\inf_{q} \left<|\pi_{b,q}|\right>.\]
Applying (\ref{est-pi-bq}) and (\ref{est-bq-l2-l3}), we obtain
\begin{equation}\notag
\begin{split}
\underline\varepsilon_b^{\frac23}\leq \left<|\pi_{b,q}|\right>^{\frac23}\lesssim &\ d_i^{\frac23}\left(\sum_{p} K_{q-p} \lambda_{p}^{2}\left<|B_{p}|^3\right>\right)^{\frac23}\\
\lesssim &\ d_i^{\frac23}\sum_{p} K_{q-p}^{\frac23} \lambda_{p}^{\frac43}\left<|B_{p}|^3\right>^{\frac23}\\
\lesssim &\ d_i^{\frac23}\sum_{p} K_{q-p}^{\frac23} \lambda_{p}^{\frac43} \lambda_p^{\frac13(3-\delta_b)} \left<|B_{p}|^2\right>\\
\lesssim &\ d_i^{\frac23}\sum_{p} K_{q-p}^{\frac23} \lambda_p^{\frac13(10-\delta_b)} \mathcal E_b(\lambda_p)\\
\end{split}
\end{equation}
where we used the scaling $\mathcal E_b(\lambda_p)\sim \frac{\left<|B_{p}|^2\right>}{\lambda_p}$. The last inequality gives an average lower bound for the energy spectrum.

\subsection{Scaling derivation of Conjecture \ref{q42}.}
\label{sec-proof-con42}

Denote the typical magnetic field difference associated with scale $\ell$ by $\delta_{\ell} B=B(x+\ell,t)-B(x,t)$. The eddy turnover time is hence 
\[t_{\ell}\sim \frac{\ell^2}{d_i\delta_{\ell} B}.\] 
Active eddies of size $\ell$ fill only a fraction $(\ell/L)^{3-\delta_b}$ of the total volume. The energy per unit mass associated with scale $\ell$ is
$E_{\ell}\sim (\delta_{\ell} B)^2 \left(\ell/L\right)^{3-\delta_b}$.
According to the energy law, we have
\begin{equation}\label{eq-sf-scaling1}
\varepsilon_b\sim E_{\ell}/t_{\ell}\sim d_i(\delta_{\ell} B)^3\ell^{-2}\left(\ell/\ell_0\right)^{3-\delta_b}.
\end{equation}
Taking $\ell=L$ in (\ref{eq-sf-scaling1}) indicates 
\begin{equation}\label{eq-sf-scaling2}
\varepsilon_b\sim d_i(\delta_{L} B)^3/L^2.
\end{equation}
Combining (\ref{eq-sf-scaling1}) and (\ref{eq-sf-scaling2}) leads to
\begin{equation}\label{eq-sf-scaling3}
\delta_{\ell} B\sim \delta_L B\left(\ell/L\right)^{(\delta_b-1)/3}.
\end{equation}
Henceforth, applying (\ref{eq-sf-scaling2}) and (\ref{eq-sf-scaling3}), the structure function is expected to satisfy
\begin{equation}\notag
\begin{split}
\left<|\delta_{\ell} B|^p\right>\sim&\ (\delta_{\ell} B)^p  \left(\ell/L\right)^{3-\delta_b}\\
\sim&\ (\delta_L B)^p\left(\ell/L\right)^{\frac{p(\delta_b-1)}{3}} \left(\ell/L\right)^{3-\delta_b}\\
\sim&\ \left(d_i^{-1}\varepsilon_{b}\right)^{\frac p3}\ell^{\frac{2p}{3}+(3-\delta_b)(1-\frac{p}{3})}
\end{split}
\end{equation}
which gives the scaling of (\ref{sf2}).


\subsection{Proof of Theorem \ref{thm-sf-bound}.}
Recall that \[\bar\varepsilon_b=\sup_{q}d_i\lambda_q^2\left<|B_q|^3\right>;\]
hence, we have 
\begin{equation}\label{est-sf-l3}
\left<|\delta_{\ell}B|^3\right>\leq d_i^{-1} \bar\varepsilon_b\ell^2.
\end{equation}
In the following, we will estimate $\left<|\delta_{\ell}B|^2\right>$ and then $\left<|\delta_{\ell}B|^p\right>$ with $2<p<3$ by interpolation.

Recall $\delta_yB(x,t)=B(x+y,t)-B(x,t)$.  Let $q$ be the integer such that $\lambda_q\sim |y|^{-1}$. We infer, by applying the Mean-value theorem
\begin{equation}\notag
\begin{split}
<|\delta_yB|^2>=&\ \frac{1}{T|\Omega|}\int_0^T\int_{\Omega}|B(x+y,t)-B(x,t)|^2\, \mathrm{d}x\mathrm{d}t\\
\leq &\ \frac{1}{T|\Omega|}\int_0^T\int_{\Omega}\sum_{p\leq q}|B_p(x+y,t)-B_p(x,t)|^2\, \mathrm{d}x\mathrm{d}t\\
&+\frac{1}{T|\Omega|}\int_0^T\int_{\Omega}\sum_{p> q}|B_p(x+y,t)-B_p(x,t)|^2\, \mathrm{d}x\mathrm{d}t\\
\lesssim &\ \frac{1}{T|\Omega|}\int_0^T\int_{\Omega}\sum_{p\leq q}|y|^2\lambda_p^2|B_p|^2\, \mathrm{d}x\mathrm{d}t+\frac{1}{T|\Omega|}\int_0^T\int_{\Omega}\sum_{p> q}|B_p|^2\, \mathrm{d}x\mathrm{d}t\\
\lesssim & \ \sum_{p\leq q}|y|^2\lambda_p^2 \left<|B_p|^2\right>+\sum_{p> q}\left<|B_p|^2\right>.
\end{split}
\end{equation}
Employing (\ref{est-bq-l2-average}) and noticing that $\ell_p\sim\lambda_p^{-1}$ and $\lambda_q\sim |y|^{-1}$, we continue with the last inequality
\begin{equation}\notag
\begin{split}
\left<|\delta_yB|^2\right>\lesssim &\ \sum_{p\leq q}|y|^2\lambda_p^2 (d_i^{-1}\bar\varepsilon_b)^{\frac23}\lambda_p^{-\frac43}\left(\frac{\ell_p}{L}\right)^{1-\frac{\delta_b}{3}}+\sum_{p> q}(d_i^{-1}\bar\varepsilon_b)^{\frac23}\lambda_p^{-\frac43}\left(\frac{\ell_p}{L}\right)^{1-\frac{\delta_b}{3}}\\
\lesssim &\ \lambda_q^{-2}\sum_{p\leq q} (d_i^{-1}\bar\varepsilon_b)^{\frac23}\lambda_p^{\frac23}\left(L\lambda_p\right)^{\frac{\delta_b}{3}-1}+\sum_{p> q}(d_i^{-1}\bar\varepsilon_b)^{\frac23}\lambda_p^{-\frac43}\left(L\lambda_p\right)^{\frac{\delta_b}{3}-1}.
\end{split}
\end{equation}
For $\delta_b\in[1,3]$, we have 
\[\frac23+\frac{\delta_b}{3}-1\geq0, \ \ \ -\frac43+\frac{\delta_b}{3}-1<0.\]
Thus, we further deduce that 
\begin{equation}\label{est-sf-bound-l2-delta}
\begin{split}
\left<|\delta_yB|^2\right>\lesssim &\ \lambda_q^{-2}(d_i^{-1}\bar\varepsilon_b)^{\frac23}\lambda_q^{\frac23}\left(L\lambda_q\right)^{\frac{\delta_b}{3}-1}+(d_i^{-1}\bar\varepsilon_b)^{\frac23}\lambda_q^{-\frac43}\left(L\lambda_q\right)^{\frac{\delta_b}{3}-1}\\
\lesssim &\ (d_i^{-1}\bar\varepsilon_b)^{\frac23}\lambda_q^{-\frac43}\left(L\lambda_q\right)^{\frac{\delta_b}{3}-1}. 
\end{split}
\end{equation}
Now, for $2<p<3$, we have by interpolation and using (\ref{est-sf-l3}) and (\ref{est-sf-bound-l2-delta})
\begin{equation}\notag
\begin{split}
\left<|\delta_yB|^p\right>\leq & \left<|\delta_yB|^2\right>^{3-p}\left<|\delta_yB|^3\right>^{p-2}\\
\lesssim &\ (d_i^{-1}\bar\varepsilon_b)^{\frac23(3-p)}\lambda_q^{-\frac43(3-p)}\left(L\lambda_q\right)^{\left(\frac{\delta_b}{3}-1\right)(3-p)}  (d_i^{-1} \bar\varepsilon_b)^{p-2}\lambda_q^{2(2-p)}\\
\lesssim &\ (d_i^{-1}\bar\varepsilon_b)^{\frac{p}3}\lambda_q^{-\frac{2p}{3}+(\delta_b-3)(1-\frac{p}{3})}\\
\lesssim &\ (d_i^{-1}\bar\varepsilon_b)^{\frac{p}3}|y|^{\frac{2p}{3}+(3-\delta_b)(1-\frac{p}{3})}.
\end{split}
\end{equation}
It completes the proof of Theorem \ref{thm-sf-bound}.

\subsection{Scaling derivation of Conjecture \ref{q43}.}
\label{sec-proof-con43}

The scalings in Conjecture \ref{q43} will be attained through heuristic analysis and estimates all the flux terms using harmonic analysis tools. 
For the magnetic field in the sub-ion range, the dynamics captured by the Hall term is active and dominant. In this scale range, the magnetic energy spectrum is expected to be similar as that of the EMHD, which is predicted in Conjecture \ref{q41}. Thus we focus on the ion-inertial range here and neglect the flux from the Hall term. 

The energy law of the Hall MHD system (\ref{hmhd}) at the $q$-th shell is given by
\begin{equation}\label{hmhd-energy-uq1}
\begin{split}
\frac12\frac{d}{dt}\int_{\mathbb T^3}|u_q|^2\, \mathrm{d}x+\int_{\mathbb T^3} (u\cdot\nabla u)_q\cdot u_q\, \mathrm{d}x-\int_{\mathbb T^3} (B\cdot\nabla B)_q\cdot u_q\, \mathrm{d}x\\
+\nu \int_{\mathbb T^3}|\nabla u_q|^2\, \mathrm{d}x=0,
\end{split}
\end{equation}
\begin{equation}\label{hmhd-energy-bq1}
\begin{split}
\frac12\frac{d}{dt}\int_{\mathbb T^3}|B_q|^2\, \mathrm{d}x+\int_{\mathbb T^3} (u\cdot\nabla B)_q\cdot B_q\, \mathrm{d}x-\int_{\mathbb T^3} (B\cdot\nabla u)_q\cdot B_q\, \mathrm{d}x\\
+d_i\int_{\mathbb T^3}\nabla\times((\nabla\times B)\times B)_q\cdot B_q\, \mathrm{d}x+\mu \int_{\mathbb T^3}|\nabla B_q|^2\, \mathrm{d}x=0.
\end{split}
\end{equation}

Denote 
\[\varepsilon_{u,q}=\nu \int_{\mathbb T^3}|\nabla u_q|^2\, \mathrm{d}x\]
 by the kinetic energy dissipation rate at the $q$-th shell; and 
 \[\varepsilon_{b,q}=\mu \int_{\mathbb T^3}|\nabla B_q|^2\, \mathrm{d}x\]
 is similarly introduced previously. 
Based on (\ref{hmhd-energy-uq1}), we have the following scaling relationships at the scale of $\lambda_q$ by considering the nonlinearity $(u\cdot\nabla)u$
\begin{equation} \label{scal-u1}
\varepsilon_{u,q}\sim  \frac12\frac{d}{dt}\|u_q\|_{L^2}^2\sim \nu\lambda_q^2\|u_q\|_{L^2}^2\sim \|u_q\nabla u_qu_q\|_{L^1}
\end{equation}
and the following scaling by considering the nonlinearity $(B\cdot\nabla)B$
\begin{equation} \label{scal-u2}
\varepsilon_{u,q}\sim \nu\lambda_q^2\|u_q\|_{L^2}^2
\sim \|B_q\nabla B_qu_q\|_{L^1}.
\end{equation}
In the sub-ion range, neglecting the Hall effect, the energy identity (\ref{hmhd-energy-bq1}) suggests
\begin{equation}\label{scal2}
\varepsilon_{b,q}\sim  \frac12\frac{d}{dt}\|B_q\|_{L^2}^2\sim \mu\lambda_q^2\|B_q\|_{L^2}^2\sim \|\nabla\times(u_q\times B_q)\cdot B_q\|_{L^1}. 
\end{equation}
It follows from (\ref{scal-u1}), H\"older's inequality and the crucial relationship (\ref{bern2}) with intermittency dependence that
\begin{equation}\notag
\varepsilon_{u,q}\sim \nu\lambda_q^2\|u_q\|_{L^2}^2\sim \lambda_q\|u_q\|_{L^2}^2\|u_q\|_{L^\infty}\sim \lambda_q^{(5-\delta_u)/2}\|u_q\|_{L^2}^3,
\end{equation}
which implies 
\begin{equation}\label{scaling-uq-l2}
\|u_q\|_{L^2}\sim \nu \lambda_q^{\frac{\delta_u-1}{2}}.
\end{equation}
Combining $\varepsilon_{u,q}\sim \nu\lambda_q^2\|u_q\|_{L^2}^2$ with (\ref{scaling-uq-l2}) also gives
\begin{equation}\label{nu-eu}
\nu\sim \varepsilon_{u,q}^{\frac13} \lambda_q^{-\frac{\delta_u+1}{3}}.
\end{equation}
Similarly, it follows from (\ref{scal-u2}), H\"older's inequality and (\ref{bern2}) that
\begin{equation}\notag
\varepsilon_{u,q}\sim \nu\lambda_q^2\|u_q\|_{L^2}^2\sim \lambda_q\|B_q\|_{L^2}^2\|u_q\|_{L^\infty}\sim \lambda_q^{\frac{5-\delta_u}{2}}\|u_q\|_{L^2}\|B_q\|_{L^2}^2
\end{equation}
and hence
\begin{equation}\label{scaling-bq-l2}
\|B_q\|_{L^2}^2\sim \nu \lambda_q^{(\delta_u-1)/2}\|u_q\|_{L^2}.
\end{equation}
Analogously, applying (\ref{scal2}), H\"older's inequality and (\ref{bern-intermit}), we have
\begin{equation}\label{scaling-ebq}
\begin{split}
\varepsilon_{b,q}\sim \mu\lambda_q^2\|B_q\|_{L^2}^2\sim&\ \lambda_q\|u_q\|_{L^2}\|B_q\|_{L^2}\|B_q\|_{L^\infty}\\
\sim&\ \lambda_q^{(5-\delta_b)/2}\|u_q\|_{L^2}\|B_q\|_{L^2}^2.
\end{split}
\end{equation}
Combining (\ref{scaling-uq-l2}), (\ref{nu-eu}), (\ref{scaling-bq-l2}) and (\ref{scaling-ebq}), we obtain 
\begin{equation}\label{mu-eub}
\mu\sim \varepsilon_u^{-\frac16} \varepsilon_b^{\frac12} \lambda_q^{-\frac{1}{4}\delta_b-\frac{1}{12}\delta_u-\frac13}.
\end{equation}
Finally, we postulate that, by employing (\ref{scaling-uq-l2}), (\ref{scaling-bq-l2}), (\ref{nu-eu}) and (\ref{mu-eub})
\begin{equation}\notag
\begin{split}
\mathcal E_u(\lambda_q)\sim&\ \|u_q\|_{L^2}^2/\lambda_q\sim \nu^2\lambda_q^{\delta_u-2}\sim (\varepsilon_{u,q})^{\frac23}\lambda_q^{\frac13(\delta_u-8)}\\
\mathcal E_b(\lambda_q)\sim&\ \|B_q\|_{L^2}^2/\lambda_q\sim \nu\mu\lambda_q^{\frac12(\delta_u+\delta_b)-2}\sim (\varepsilon_{u,q})^{\frac16}(\varepsilon_{b,q})^{\frac12}\lambda_q^{\frac1{12}(\delta_u+3\delta_b)-\frac83}
\end{split}
\end{equation}
which inspire the predictions in Conjecture \ref{q43} by noticing both $\lambda_q$ and $k$ are wavenumber notations.

\subsection{Scaling derivation of Conjecture \ref{q44}.}
\label{sec-proof-con44}

It follows from (\ref{nu-eu}) that
\begin{equation}\notag
\lambda_q\sim (\nu^{-3}\varepsilon_{u,q})^{\frac{1}{\delta_u+1}}
\end{equation}
which inspires the scaling (\ref{hmhd-kd-u}). While, (\ref{mu-eub}) implies
\begin{equation}\notag
\lambda_q\sim \left(\mu^{-3} \varepsilon_{u,q}^{-\frac12}\varepsilon_{b,q}^{\frac32}\right)^{\frac{1}{\frac{3\delta_b+\delta_u}{4}+1}}
\end{equation}
which is the reason of the prediction of (\ref{hmhd-kd-b}).
Physical evidence shows that $\kappa_{\mathrm i}^b$ should have the scale of $1/d_i$. 






\section{Phenomenologies of intermittent turbulence for classical MHD}
\label{sec-mhd}

Since the establishment of Iroshnikow-Kraichnan (IK) \cite{Ir} scaling theory for MHD turbulence in 1960s based on K41 theory, outstanding breakthrough has been made in this area during the last few decades. Goldreich and Sridhar (GS) \cite{GS1} developed a theory with the awareness of the anisotropic feature of MHD turbulent motions at small scales and based on the so-called critical balance conjecture in 1990s. Later on, Boldyrev's theory \cite{Bold1} based on dynamic alignment in 2000s turns out to be appealing and receives lots of attention. However, Beresnyak \cite{Ber1} brought up the objection that Boldyrev's alignment theory violates the scaling symmetry of the system and fails at small scales in early 2010s. Up to now, there are still serious debates over which theory is more feasible. Nevertheless, Boldyrev's theory can be viewed as a qualitative theory of intermittency for the turbulent MHD system. In this part of the project, we derive some phenomenological laws for the MHD turbulence with dependence on the intermittency dimension.

\subsection{Energy spectrum and structure functions for MHD turbulence with intermittency correction}
\label{sec-mhd1}
The incompressible MHD system (\ref{hmhd}) with $d_i=0$ can be formulated into a more symmetric form by using Els\"asser variables 
\[Z^+=u+B, \ \ Z^-=u-B\]
introduced by Els\"asser \cite{El}. Let $B_0$ be the background mean magnetic field. Let $\rho_0$ and $\mu_0$ be the constant density of plasma and the permeability of free space. The Alfv\'en speed is given by $v_A=B_0/(\rho_0\mu_0)^{1/2}$.
The Els\"asser variables satisfy the system
\begin{equation}\label{El}
\begin{split}
\partial_tZ^\pm \mp(B_0\cdot\nabla)Z^\pm+(Z^\mp\cdot\nabla) Z^\pm+\nabla p=&\ \eta^+\Delta Z^\pm+\eta^-\Delta Z^\mp, \\
 \nabla\cdot Z^\pm=&\ 0, 
\end{split}
\end{equation}
with $\eta^\pm=(\nu\pm\mu)/2$. 
In principle, the Els\"asser variables represent Alfv\'en wave perturbations propagating in the direction or the opposite direction of the background magnetic field $B_0$. 
 There are evidences both from theoretical analysis \cite{DMV} and in-situ data observations \cite{BCa} that it is useful to apply the Els\"asser formulation (\ref{El}) in the study of plasma turbulence.
 
 \begin{Remark}
Since we are interested in plasmas with high magnetic Reynolds number which indicates $\mu$ is small enough, it is thus natural to assume $\eta^->0$.
\end{Remark}

 

In the means of Definition \ref{int-def}, we define the intermittency dimension $\delta^+$ for the field $Z^+$, and the intermittency dimension $\delta^-$ for $Z^-$. In 3D, we have $\delta^\pm\in[0,3]$ and the crucial relations at the scale $\lambda_q$
\begin{equation}\label{bern-z}
\|Z^+_q\|_{L^\infty}\sim \lambda_q^{(3-\delta^+)/2} \|Z^+_q\|_{L^2}, \ \ \|Z^-_q\|_{L^\infty}\sim \lambda_q^{(3-\delta^-)/2} \|Z^-_q\|_{L^2}.
\end{equation}
Denote the energy quantities associated with each Els\"asser variable as
\[E^\pm(t)=\frac12\|Z^\pm(t)\|_{L^2}^2.\] 
Denote
 \[\varepsilon^\pm=\eta^\pm\left<|\nabla Z^\pm|^2\right>\]
as the average energy dissipation rate of $Z^\pm$ and 
$\mathcal E^\pm(k)$ as the energy spectrum.  
We have the following predictions on the energy spectra scaling for the turbulent fields $Z^\pm$.

\begin{Conjecture}\label{q51}
{\it The energy spectra $\mathcal E^\pm(k)$ in the inertial range for the 3D MHD model (\ref{El}) in Els\"asser formulation obey the scaling law }
\begin{equation}\label{en-in}
\begin{split}
\mathcal E^+(k)\sim \left((\ea)^2/\eb\right)^{2/3}k^{(2\deb-\dea-8)/3}, \\
\mathcal E^-(k)\sim \left((\eb)^2/\ea\right)^{2/3}k^{(2\dea-\deb-8)/3}.
\end{split}
\end{equation}
\end{Conjecture}

When both fields $Z^\pm$ are homogeneous isotropic and self-similar,  i.e. $\delta^\pm=3$, the scaling law becomes
$\mathcal E^\pm(k)\sim \left((\varepsilon^\pm)^2/\varepsilon^\mp\right)^{2/3}k^{-5/3}$.
It is worth to notice that this scaling is consistent with Goldreich-Sridhar theory \cite{GS1} in the perpendicular cascade.
In the extreme anisotropic case of $\delta^\pm=0$, the scaling (\ref{en-in}) becomes $\mathcal E^\pm(k)\sim \left((\varepsilon^\pm)^2/\varepsilon^\mp\right)^{2/3}k^{-8/3}$.
Therefore, the conjecture gives a good explanation why experimental data showed energy spectra between $k^{-5/3}$ and $k^{-8/3}$ for different turbulent fields. The crucial point is that turbulent fields may have different intermittency levels.

Upper and lower bounds for $\mathcal E^\pm(k)$ are expected to be attained in an analogous way as for Theorem \ref{thm-bounds-spectrum}. 

We explain below how to attain the prediction (\ref{en-in}) by using scaling analysis and harmonic analysis techniques in analogy with the analysis of Conjecture \ref{q43}. 
The energy law of (\ref{El}) 
implies the following scaling relations at each scale, 
\[\varepsilon^\pm\sim \eta^\pm\|\nabla Z_q^\pm\|_{L^2}^2\sim \|(Z_q^\mp\cdot \nabla) Z_q^\pm\cdot Z_q^\pm\|_{L^1}.\]
Based on Littlewood-Paley theory in harmonic analysis, it is trivial to see 
\begin{equation}\notag
\eta^\pm\|\nabla Z_q^\pm\|_{L^2}^2\sim \eta^\pm\lambda_q^2\|Z_q^\pm\|_{L^2}^2.
\end{equation}
On the other hand, thanks to (\ref{bern-z}), it follows from H\"older's inequality and Littlewood-Paley theory that
\[\|(Z_q^\mp\cdot \nabla) Z_q^\pm\cdot Z_q^\pm\|_{L^1}\sim 
\lambda_q\|Z_q^\pm\|_{L^2}^2\lambda_q^{(3-\delta^\mp)/2}\|Z_q^\mp\|_{L^2}. \]
Combining the above relations leads to
\[\|Z_q^\pm\|_{L^2}\sim \eta^\mp\lambda_q^{(\delta^\pm-1)/2}, \ \ \
\varepsilon^\pm\sim \eta^\pm(\eta^\mp)^2\lambda_q^{1+\delta^\pm}. \]
Thus, we infer by scaling that 
\[\mathcal E^\pm(\lambda_q)\sim \frac{\|Z^\pm_q\|_{L^2}^2}{\lambda_q} 
\sim \left((\varepsilon^\pm)^2/\varepsilon^\mp\right)^{\frac23}\lambda_q^{(2\delta^\mp-\delta^\pm-8)/3} \ \ \Longrightarrow \ \ \mbox {inspires Conjecture \ref{q51}}. \]

Further heuristic scaling analysis leads to the prediction of the cutoff scaling between the dissipation range and the inertial range. 

\begin{Conjecture}\label{q52}
{\it There exists a critical dissipation wavenumber $\kappa_{\mathrm d}^+$ for the field $Z^+$ and a critical dissipation wavenumber $\kappa_{\mathrm d}^-$ for $Z^-$ such that $\kappa_{\mathrm d}^{\pm}$ separate the dissipation regime from the inertial regime for $Z^\pm$ respectively. 
Moreover, we predict that }
\begin{equation}\label{kd}
\begin{split}
\kappa_{\mathrm d}^+\sim \left(\ea/[(\nu+\mu)(\nu-\mu)^2]\right)^{1/(1+\dea)}, \\
\kappa_{\mathrm d}^-\sim \left(\eb/[(\nu-\mu)(\nu+\mu)^2]\right)^{1/(1+\deb)}.
\end{split}
\end{equation}
\end{Conjecture}

We notice the dissipation wavenumber $\kappa_{\mathrm d}^+$ depends on the intermittency dimension $\dea$ of $Z^+$, not on $\deb$; vice versa for $\kappa_{\mathrm d}^-$. Nevertheless, both $\kappa_{\mathrm d}^\pm$ depend on both of the fluid and magnetic Reynolds number $1/\nu$ and $1/\mu$.
The motivation of (\ref{kd}) comes from the early derived relation $\varepsilon^\pm\sim \eta^\pm(\eta^\mp)^2\lambda_q^{1+\delta^\pm}$ which is equivalent to 
$\lambda_q\sim \left(\varepsilon^\pm/\eta^\pm(\eta^\mp)^2\right)^{1/(1+\delta^\pm)}$. 


Regarding structure functions, we have:

\begin{Conjecture}\label{q53}
{\it
The $p$-th order structure functions $\left<|\delta Z^\pm_{\ell}|^p\right>$ for intermittent fields $Z^\pm$ obeys the scaling law
\begin{equation}\label{sf-Z-law}
\left<|\delta_{\ell} Z^\pm|^p\right>
\sim \left(\varepsilon^\pm\right)^{\frac p3}\ell^{\frac{p}{3}+(3-\delta^\pm)(1-\frac{p}{3})}.
\end{equation}
}
\end{Conjecture}
We present a brief heuristic analysis below to produce (\ref{sf-Z-law}).
Denote
\[\delta_{\ell} Z^\pm=Z^\pm(x+\ell,t)-Z^\pm(x,t)\] 
by the difference of $Z^\pm$ associated with scale $\ell$.
The eddy turnover time is hence 
\[t_{\ell}\sim \frac{\ell}{\delta_{\ell} Z^\pm}.\] 
For $Z^\pm$, active eddies of size $\ell$ fill only a fraction $(\ell/L)^{3-\delta^\pm}$ of the total volume. The energy per unit mass associated with scale $\ell$ is
$E^\pm_{\ell}\sim (\delta_{\ell} Z^\pm)^2 \left(\ell/L\right)^{3-\delta^\pm}$.
According to the energy law, we have
\begin{equation}\label{eq-sf-Z-scaling1}
\varepsilon^\pm\sim E^\pm_{\ell}/t_{\ell}\sim (\delta_{\ell} Z^\pm)^3\ell^{-1}\left(\ell/L\right)^{3-\delta^\pm}.
\end{equation}
Taking $\ell=L$ in (\ref{eq-sf-Z-scaling1}) indicates 
\begin{equation}\label{eq-sf-Z-scaling2}
\varepsilon^\pm\sim (\delta_{L} Z^\pm)^3/L.
\end{equation}
Combining (\ref{eq-sf-Z-scaling1}) and (\ref{eq-sf-Z-scaling2}) leads to
\begin{equation}\label{eq-sf-Z-scaling3}
\delta_{\ell} Z^\pm\sim \delta_L Z^\pm\left(\ell/L\right)^{(\delta^\pm-2)/3}.
\end{equation}
Henceforth, applying (\ref{eq-sf-Z-scaling2}) and (\ref{eq-sf-Z-scaling3}), the structure function is expected to satisfy
\begin{equation}\notag
\begin{split}
\left<|\delta_{\ell} Z^\pm|^p\right>\sim&\ (\delta_{\ell} Z^\pm)^p  \left(\ell/L\right)^{3-\delta^\pm}\\
\sim&\ (\delta_L Z^\pm)^p\left(\ell/L\right)^{\frac{p(\delta^\pm-2)}{3}} \left(\ell/L\right)^{3-\delta^\pm}\\
\sim&\ \left(\varepsilon^\pm\right)^{\frac p3}\ell^{\frac{p}{3}+(3-\delta^\pm)(1-\frac{p}{3})}
\end{split}
\end{equation}
which gives the scaling of (\ref{sf-Z-law}).



\subsection{Perpendicular cascade of intermittent MHD}
The assumption that MHD turbulence consists of perturbations with $k_{\perp}\gg k_{||}$ but the Alfv\'enic propagation remains important leads to the equations:
\begin{equation}\label{eq-El-2}
\begin{split}
\partial_tZ^\pm_{\perp}\mp v_A\nabla_{||}Z^\pm_{\perp}+Z^\mp_{\perp}\cdot\nabla_{\perp}Z^\pm_{\perp}+\nabla_{\perp}p=&\ \eta^\pm\Delta_{\perp}Z^\pm_{\perp}, \\
\nabla_{\perp}\cdot Z^\pm_{\perp}=&\ 0.
\end{split}
\end{equation}
Since $Z^\pm_{\perp}$ oscillate in $\mathbb R^2$, their intermittency dimensions $\delta_{\perp}^\pm$ can be defined similarly as Definition \ref{int-def}  with dimension $n=2$ and hence $\delta_{\perp}^\pm\in [0,2]$. The $L^\infty$ norm and $L^2$ norm are thus related at each scale as
\begin{equation}\label{bern-zper}
\|Z^+_{q\perp}\|_{L^\infty}\sim \lambda_q^{(2-\delta_{\perp}^+)/2} \|Z^+_{q\perp}\|_{L^2}, \ \ \|Z^-_q\|_{L^\infty}\sim \lambda_{q\perp}^{(2-\delta_{\perp}^-)/2} \|Z^-_{q\perp}\|_{L^2}.
\end{equation}

Denote the energy quantities associated with each Els\"asser variable perpendicular to the background magnetic field by
$E_{\perp}^\pm(t)=\frac12\|Z_{\perp}^\pm(t)\|_{L^2}^2$. 
Let $\varepsilon_{\perp}^\pm$  stand  for the average energy dissipation rate of $E_{\perp}^\pm(t)$. Let $\mathcal E^\pm(k_{\perp})$ represent the energy spectrum corresponding to $Z^\pm$. 

Analogous heuristic analysis and energy law of $Z^\pm$ lead to the following scaling law. 

\begin{Conjecture}\label{q54}
{\it
Assume the nonlinear interactions dominate over the Alfv\'enic propagation.
 Then the energy spectra $\mathcal E^\pm(k_{\perp})$ in their inertial range for the 3D MHD model (\ref{eq-El-2}) obey the scaling law 
\begin{equation}\notag
\begin{split}
\mathcal E^+(k_{\perp})\sim \left((\varepsilon_{\perp}^+)^2/\varepsilon_{\perp}^-\right)^{\frac23}k_{\perp}^{(2\delta_{\perp}^--\delta_{\perp}^+-7)/3}, \\
\mathcal E^-(k_{\perp})\sim \left((\varepsilon_{\perp}^-)^2/\varepsilon_{\perp}^+\right)^{\frac23}k_{\perp}^{(2\delta_{\perp}^+-\delta_{\perp}^--7)/3}.
\end{split}
\end{equation}
}
\end{Conjecture}

It is important to notice that for homogeneous isotropic and self-similar turbulent fields, i.e., $\delta_{\perp}^\pm=2$, the scaling exponent is $-5/3$ which recovers GS scaling \cite{GS1}, see Figure \ref{fig5}. In general, assume $\delta_{\perp}^+=\delta_{\perp}^-=\delta_{\perp}$, then the scaling is $\mathcal E^\pm(k_{\perp})\sim k^{-\gamma}_{\perp}$ with $\gamma=(7-\delta_{\perp})/3$, see Figure \ref{fig6}.

\begin{minipage}{.48\textwidth}
\begin{center}
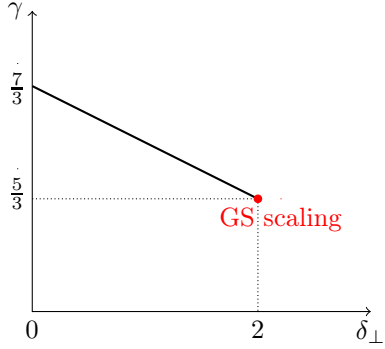

\begin{tikzpicture}
\draw [<->] (0,4) node[left]{$\gamma$}  --(0,0)  -- (4.5,0) node[below]{$\delta_{\perp}$};
\draw  [thick] (0,3)  -- (3,1.5) ; 
\draw [densely dotted] (3, 1.5)  -- (3,0) ; 
\draw [densely dotted] (0, 1.5)  -- (3,1.5) ;  
\draw  (-0.2, 1.9) -- node[below]{$\frac53$}  (-0.2, 1.9)  ;
\draw  (0,0) -- node[below]{$0$}  (0,0)  ;
\draw  (3,0) -- node[below]{$2$}  (3,0)  ;
\draw  [red] (3.3,1.5) -- node[below]{GS scaling}  (3.3,1.5)  ;
\draw  (-0.2, 3.3) -- node[below]{$\frac73$}  (-0.2, 3.3)  ;
\draw [fill, red] (3,1.5) circle [radius=0.05];
\end{tikzpicture}
\captionof{figure}{Negative exponent $\gamma$ of perpendicular energy spectrum with dependence on intermittency dimension $\delta_{\perp}$.}
\label{fig5}
\end{center}
\end{minipage}
\begin{minipage}{.48\textwidth}
\begin{center}
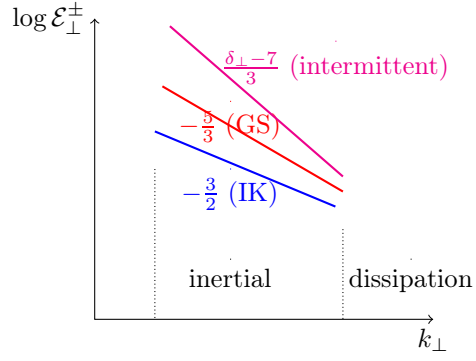

\begin{tikzpicture}
\draw [<->] (0,4) node[left]{$\log\mathcal E_{\perp}^\pm$}  --(0,0)  -- (4.5,0) node[below]{$k_{\perp}$};
\draw  [thick, blue] (0.8,2.5)  -- (3.2,1.5) ;  
\draw  [blue] (1.8,2) -- node[below]{$-\frac32$ (IK)}  (1.8,2)  ;
\draw  [thick, red] (0.9,3.1)  -- (3.3,1.7) ; 
\draw  [red] (1.8,2.9) -- node[below]{$-\frac53$ (GS)}  (1.8,2.9)  ;
\draw  [thick, magenta] (1,3.9)  -- (3.3,1.9) ; 
\draw [magenta] (3.2, 3.7) -- node[below]{$\frac{\delta_{\perp}-7}{3}$ (intermittent)}  (3.2, 3.7)  ;
\draw [densely dotted] (0.8, 2)  -- (0.8,0) ; \draw [densely dotted] (3.3, 1.2)  -- (3.3,0) ; 
\draw  (1.8,0.8) -- node[below]{inertial}  (1.8,0.8)  ;
\draw  (4.2,0.8) -- node[below]{dissipation}  (4.2,0.8)  ;
\end{tikzpicture}
\captionof{figure}{Energy spectra of perpendicular cascade under different assumptions.}
\label{fig6}
\end{center}
\end{minipage}

\begin{Conjecture}\label{q55}
{\it
There exist critical dissipation wavenumber $\kappa_{\mathrm d {\perp}}^\pm$ for the fields $Z^\pm$
such that $\kappa_{\mathrm d \perp}^{\pm}$ separate the dissipation regime from the inertial regime for $Z^\pm$. Moreover, we predict that
\begin{equation}\label{kdper}
\begin{split}
\kappa_{\mathrm d\perp}^+\sim \left(\varepsilon_{\perp}^+/[(\nu+\mu)(\nu-\mu)^2]\right)^{1/(2+\delta_{\perp}^+)}, \\
\kappa_{\mathrm d\perp}^-\sim \left(\varepsilon_{\perp}^-/[(\nu-\mu)(\nu+\mu)^2]\right)^{1/(2+\delta_{\perp}^-)}.
\end{split}
\end{equation}
}
\end{Conjecture}

We notice that $\kappa_{\mathrm d\perp}^+$ depends on $\delta_{\perp}^+$ not on $\delta_{\perp}^-$; vice versa for $\kappa_{\mathrm d\perp}^-$. Nevertheless, both $\kappa_{\mathrm d\perp}^\pm$ depend on both of the velocity and magnetic Reynolds number. 

The next conjecture recovers Boldyrev's (and IK) scaling as well as GS scaling (see Figure \ref{fig6}).

\begin{Conjecture}\label{q56}
{\it
Assume the Alfv\'enic propagation is strong compared to the nonlinear interactions. The energy spectra in the inertial regime satisfy
\begin{equation}\notag
\mathcal E^\pm(k_{\perp})\sim
\begin{cases}
 (\varepsilon_{\perp}^\pm)^{\frac23}k_{\perp}^{-\frac53}, \ \ \ \ \ \ \ \ \mbox{if the critical balance conjecture} \ \ \ell_{\perp}\sim \varepsilon ^{\frac12}\left(v_A/\ell_{||}\right)^{-\frac32}\\
\ \ \ \ \ \ \ \ \ \ \ \ \ \ \ \ \ \ \ \ \ \ \mbox {holds,}\\ 
 (\varepsilon_{\perp}^\pm v_A)^{\frac12}k_{\perp}^{-\frac32}, \ \ \ \ \mbox{if the alignment condition}\ \ \ell_{||}/\ell_{\perp}\sim v_A^2/Z_{\ell_\perp}^2 \ \ \mbox{holds.}
\end{cases}
\end{equation}
}
\end{Conjecture}
Remarkably, the second case recovers IK and Boldyrev's scaling. While under the assumption of critical balance conjecture, it is consistent with GS scaling and Conjecture \ref{q54} with $\delta_{\perp}^\pm=2$.

To support Conjecture \ref{q56}, we extract the Alfv\'enic propagation part and dissipation part from the energy law and deduce
\[\varepsilon_{\perp}^\pm \sim \eta^\pm\|\nabla_{\perp} Z_{q\perp}^\pm\|_{L^2}^2\sim v_A\|\nabla_{||} Z_{q\perp}^\pm\cdot Z_{q\perp}^\pm\|_{L^1}.\]
Performing scaling analysis and using harmonic analysis tools, we infer 
\begin{equation}\label{zper-relation2}
v_A\lambda_{q||}\sim \eta^+\lambda_{q\perp}^2\sim \eta^-\lambda_{q\perp}^2, \ \ \
\|Z_{q\perp}^\pm\|_{L^2}^2\sim (\eta^{\pm})^{-1}\varepsilon_{\perp}^\pm\lambda_{q\perp}^{-2}.
\end{equation}
It is worth to point out that (\ref{zper-relation2}) shows the relationship between the changes of the parallel and perpendicular scales. 
The critical balance conjecture along with (\ref{zper-relation2}) gives rise to 
\[\mathcal E^\pm(\lambda_{q\perp})\sim \|Z^\pm_{q\perp}\|_{L^2}^2/\lambda_{q\perp}
\sim (\varepsilon_{\perp}^\pm)^{2/3}\lambda_{q\perp}^{-5/3}.\]
On the other hand, the alignment condition $v_A^2/Z_{\ell_\perp}^2\sim\ell_{||}/\ell_{\perp}\sim \lambda_{q\perp}/\lambda_{q||}$ 
together with (\ref{zper-relation2}) implies
\[\mathcal E^\pm(\lambda_{q\perp})\sim \|Z^\pm_{q\perp}\|_{L^2}^2/\lambda_{q\perp} \sim (\varepsilon_{\perp}^\pm v_A)^{1/2}\lambda_{q\perp}^{-3/2}.\]

\end{document}